\documentclass[11pt]{article}
\usepackage{amssymb}
\usepackage{epsfig}
\usepackage{mathrsfs}

\def\eps{\varepsilon}
\def\ZZ{\mathbb Z}
\def\RR{\mathbb R}

\def\na{\nabla}
\def\pa{\partial}

\def\ds{\displaystyle}
\def\vphi{\varphi}

\def\Se{{\Sigma_{\eps}}}
\def\wSe{{\widetilde \Sigma_{\eps}}}
\def\wTe{{\widetilde T_{\eps}}}
\def\tgamma{\widetilde\gamma}
\def\Rn{\RR^{n+1}_+}
\def\tw{\widetilde w}

\def\Wa{W^{1,2}(D,|y|^a)}
\def\Wal{W_{\rm{loc}}^{1,2}(D,|y|^a)}
\def\Wao{W_0^{1,2}(D,|y|^a)}

\def\div{\mbox{div}\,}

\def\ov{\overline v}
\def\ou{\overline u}

\def\ovae{{\overline v}_{\alpha}^\eps}
\def\vae{{ v}_{\alpha}^\eps}

\def\wae{\overline {w}_{\alpha}^\eps}

\def\om{\overline m}

\def\oell{\overline \ell}

\def\cor{ w^\eps}

\def\J{\mathscr J}
\def\Ja{\mathscr J_\alpha}
\def\capa{{\rm cap}}

\newtheorem{theorem}{Theorem}[section]
\newtheorem{lemma}[theorem]{Lemma}
\newtheorem{proposition}[theorem]{Proposition}

\newtheorem{remark}[theorem]{Remark}
\newtheorem{corollary}[theorem]{Corollary}

\def\qed{\hbox{${\vcenter{\vbox{                        
   \hrule height 0.4pt\hbox{\vrule width 0.4pt height 6pt
   \kern5pt\vrule width 0.4pt}\hrule height 0.4pt}}}$}}

\title{Random Homogenization of Fractional Obstacle Problems}

\author{L. A. Caffarelli\footnote{Dept of Mathematics, University of Texas at Austin,  Austin, TX~78712, USA} \quad
and A. Mellet\footnote{Dept of Mathematics, University of British Columbia, Vancouver, BC V6T 1Z2, Canada}}

\date{}

\begin{document}

\maketitle

 \begin{abstract}
We use a characterization of the fractional Laplacian as a Dirichlet to Neumann operator 
for an appropriate differential equation to study its obstacle problem in perforated domains.
 \end{abstract}

\section{Introduction}
Given a smooth function  $\vphi:\RR^n\mapsto \RR^n$ and 
a subset  $T_\eps$ of $\RR^n$, we consider $v^\eps(x)$ solution of the following obstacle problem:
\begin{equation} \label{eq:1}
\left\{
\begin{array}{ll}
v^\eps(x) \geq \vphi(x) & \mbox{ for } x\in T_\eps  \\[5pt]
(-\Delta)^{s} v^\eps \geq 0 &  \mbox{ for } x\in \RR^n  \\[5pt]
(-\Delta)^{s} v^\eps = 0 & \mbox{ for }x\in \RR^n\setminus T_\eps \mbox{ and for }x \in T_\eps \mbox{ if } v^\eps(x) > \vphi(x).
\end{array}
\right.
\end{equation}
The operator $(-\Delta)^{s}$ denotes the fractional Laplace operator of order $s$, where $s$ is a real number between $0$ and $1$. It can be defined using Fourier transform, by  $\mathcal F((-\Delta)^{s}f)(\xi) = |\xi|^{2s} \widehat f(\xi)$. In particular, (\ref{eq:1}) can be seen as the Euler-Lagrange equation for the minimization of the $\stackrel{.}{H}^s$ norm $||f || _{\stackrel{.}{H}^s} = || \widehat f (\xi)|\xi|^s ||_{L^2}$ with the constrain that $f\geq \vphi$ on $T_\eps$.
We will see that this system of equations can also be stated 
as a boundary obstacle problem for elliptic degenerate equations.

In (\ref{eq:1}), the domain $\RR^n$ is perforated and the obstacle $\vphi(x)$ is viewed by $v^\eps(x)$ only on the subset $T_\eps$. 
A typical example of $T_\eps$ is given by:
\begin{equation}\label{eq:Te}
T_\eps = \bigcup_{k\in\ZZ^n} B_{a^\eps}(\eps k),
\end{equation}
with $a^\eps \ll \eps$.
The goal of this paper  is  to study  the asymptotic behavior of $v^\eps$ as $\eps \rightarrow 0$.
When $T_\eps$ is given by (\ref{eq:Te}),
the effective equation satisfied by the limit of  $v^\eps$ strongly  depends on the radius $a^\eps$:
If  $a^\eps$ is large enough, the limit turns out to be an obstacle problem with  obstacle $\vphi(x)$. On the other hand, if $a^\eps$ is small then the limiting problem is a simple elliptic equation without any obstacle condition. 
It is well known in the case of the regular Laplace operator ($s=1$) that there is a critical size for $a^\eps$ for which interesting behavior arises.

In the case of the regular Laplace operator,
this problem was first studied for  periodic $T_\eps$ by L. Carbone and F. Colombini \cite{CC} and then in a more general framework by E. De Giorgi, G. Dal Maso and P. Longo \cite{DDL} and G. Dal Maso and P. Longo \cite{DL}, G. Dal Maso \cite{D}.
Our main reference will be the papers of D. Cioranescu and F. Murat \cite{CM1,CM2}, in which the case of  a periodic distribution of balls  is studied.
More precisely, they prove that when $s=1 $ and if $T_\eps$ is given by  (\ref{eq:Te}) with
$a^\eps = r_0 \,\eps^\frac{n}{n-2}$, then the function $v=\lim_{\eps\rightarrow 0} v^\eps$ solves
$$ -\Delta v - \mu (v-\vphi)_-=0,$$
where $\mu$ is a real number (depending on $r_0$) and $w_-=\max(-w,0)$.
The obstacle condition thus disappears when $\eps$ goes to zero, but it gives rise to a new term $\mu (v-\vphi)_-$ in the equation.

\vspace{10pt}

In \cite{CM}, we generalize this result (still with  $s=1$) to sets $T_\eps$ that are the union of small sets $S_\eps(k)\subset\RR^n$ still periodically distributed, but with random sizes and shapes.
 More precisely, 
we introduce a probability space $(\Omega,\mathcal F, \mathcal P)$ and we 
assume that for  every $\omega\in\Omega$ and every $\eps>0$ we are given some subsets $S_\eps(k,\omega)$ such that
$$S_\eps(k,\omega) \subset B_\eps(\eps k).$$
We then  consider
$$T_\eps(\omega) = \bigcup_{k\in\ZZ^n} S_\eps(k,\omega) .$$
The only assumptions necessary to generalize the result of D. Cioranescu and F. Murat \cite{CM1}-\cite{CM2} is that each set $S_\eps(k,\omega) $ is of capacity of order $\eps^n$: $\capa (S_\eps(k,\omega))=\eps^n \gamma(k,\omega)$ (this is where the critical exponent $\eps^\frac{n}{n-2}$ comes from) and that
the $\gamma(k,\omega)$ have some averaging properties (stationary ergodicity).
\vspace{10pt}

In the present paper, we  extend the result of \cite{CM} to the case of fractional Laplace operators $s\in(0,1)$.
We will show that under appropriate assumptions on the size of the sets $S_\eps(k,\omega)$, the function  $v(x)=\lim_{\eps\rightarrow 0} v^\eps(x)$ solves
$$ (-\Delta)^{s} v - \mu (v-\vphi)_-=0.$$
In the particular  case of sets $T_\eps$ of the form (\ref{eq:Te}), the critical size is now given by 
$$a^\eps =   r_0\, \eps^\frac{n}{n-2s}$$
(the critical exponent $\frac{n}{n-2s}$  is related to the $s$-capacity of the sets $S_\eps(k,\omega)$).
\vspace{10pt}

In the remainder of this section, we briefly motivate the problem and we introduce the extension problem for the fractional Laplace operators, which allows us to rewrite  (\ref{eq:1}) as a boundary obstacle problem for a local (degenerate) elliptic operator.
The precise hypothesis on $T_\eps(\omega)$ will be detailed in the following section in which the precise statement of the main theorem is also given. The remainder of the paper is devoted to the proof of our main statement.

\subsection{A semipermeable membrane problem.}
When $s=1/2$, (\ref{eq:1}) naturally arises as a boundary  obstacle problem for the regular Laplace operator (also know as Signorini problem):
We consider the following problem set in the upper-half space $ \RR^{n+1}_+ = \{(x,y)\in \RR^n\times \RR\, ;\, y \geq 0\}$:
\begin{equation}\label{eq:2}
\left\{
\begin{array}{ll}
 - \Delta u(x,y) = 0 & \mbox{ for } (x,y)\in \RR^{n+1}_+  \\[5pt]
 u (x,0) \geq \vphi(x)  & \mbox{ for } x\in T_\eps \\[5pt]
 \pa_y u (x,0)\;\leq \;0 & \mbox{ for } x\in \RR^n \\[5pt]
 \pa_y u (x,0)\;= \; 0 & \mbox{ for } x\in \RR^n \setminus T_\eps \mbox{  and for }  x\in T_\eps \mbox{ if } u(x,0)>\vphi(x) 
\end{array}
\right.
\end{equation}
with the boundary condition 
$$\lim_{y\rightarrow \infty} u(x,y) =0.$$
It is then well-known that  $v(x)=u(x,0)$  is solution of (\ref{eq:1}) with $s=1/2$ (see \cite{S} and \cite{CSS} for details).
\vspace{10pt}

It can be of interest to state equation (\ref{eq:2}) in a bounded domain $D\subset \RR^{n+1}_+$:
Introducing
$$\Sigma = D \cap \{y=0\}\quad \mbox{ and } \quad\Gamma = \pa D \cap\{y>0\},$$
we can consider the following boundary obstacle problem:
\begin{equation}\label{eq:2bis}
\left\{
\begin{array}{ll}
 - \Delta u(x,y) = 0 & \mbox{ for } (x,y)\in D  \\[5pt]
 u (x,0) \geq \vphi(x)  & \mbox{ for } x\in \Sigma\cap T_\eps \\[5pt]
 \pa_y u (x,0)\;\leq \;0 & \mbox{ for } x\in \Sigma\\[5pt]
 \pa_y u (x,0)\;= \; 0 & \mbox{ for } x\in \Sigma\setminus T_\eps \mbox{  and for }  x\in T_\eps \mbox{ if } u(x,0)>\vphi(x) 
\end{array}
\right.
\end{equation}
with the boundary condition 
$$u(x,y)= g(x,y)  \mbox{ for } (x,y)\in   \Gamma.$$

Equation (\ref{eq:2bis}) arises, for instance, in the modeling of diffusion through semi-permeable membranes (such as the membrane of a cell): The membrane is modeled by the surface $\Sigma$. The outside concentration of molecules is given by $\vphi(x)$, and the transport of molecules through the membrane and in the direction of the concentration gradient is possible  only across some given channels (represented by the set $T_\eps$) and only from the outside of the cell ($\{y<0\}$) toward the inside of the cell $D$. 
At equilibrium, the concentration inside the cell is then given by the solution $u(x,y)$ of (\ref{eq:2bis}).

\subsection{An extension problem for fractional obstacle problems}
Following L. Caffarelli, S. Salsa and L. Silvestre \cite{CSS}, we can actually rewrite (\ref{eq:1}) as a boundary obstacle problem for all fractional powers $s\in(0,1)$.
We rely for this  on the following  extension formula established by L.~Caffarelli and L. Silvestre \cite{CS}: 
For a given function $f(x)$ defined in $\RR^n$, if we define $u(x,y)$ by
\begin{equation}\label{eq:extension}
\left\{
\begin{array}{ll}
 - \div(y^a \na u) = 0 & \mbox{ for } (x,y)\in \RR^{n+1}_+ \\[5pt]
 u(x,0)= f (x)& \mbox{ for }x\in  \RR^n,
\end{array}
\right.
\end{equation}
then
$$
(-\Delta)^s f (x)= \lim_{y\rightarrow 0} y^a \pa_y u (x,y)
$$
with  
$$s=(1-a)/2.$$

We can thus rewrite the fractional obstacle problem (\ref{eq:1}) as follows:
\begin{equation} \label{eq:3}
\left\{
\begin{array}{ll}
 - \div( y^a \na  u^{\eps}) = 0 & \mbox{ for }(x,y)\in \RR^{n+1}_+ \\[5pt]
u^{\eps}(x,0) \geq \vphi (x) & \mbox{ for }x\in T_{\eps} \\[5pt]
\ds \lim_{y\rightarrow 0 }y^a  \pa_y u^{\eps}  (x,y)\leq 0 & \mbox{ for } x\in \RR^n \\[8pt]
\ds \lim_{y\rightarrow 0 }y^a  \pa_y u^{\eps}  (x,y)=  0 & \mbox{ for } x \in  \RR^n\setminus T_\eps \mbox{ and } x \in   T_\eps \cap\{u^\eps > \vphi\}
\end{array}
\right.
\end{equation}
where $a=1-2s$ (note that $a\in(-1,1)$).
Our main result will concern problems such as (\ref{eq:3}) with possibly bounded domain $D$ instead of $\RR^{n+1}_+$.

In the sequel, the theory of degenerate elliptic equations in weighted Sobolev spaces will play an important role. We refer to \cite{FKS} for many results that will be used in this paper.

\subsection{Variational formulation}
The system of equations (\ref{eq:3}) can also be written as a minimization problem.
For a given open subset $D$ of $\RR^{n+1}_+$, we denote by $L^2(D,|y|^a)$ the weighted $L^2$ space with weight $|y|^a$ and  by $W^{1,2}(D,|y|^a)$ the corresponding Sobolev's space.
We have
$$ ||u||^2_{W^{1,2}(D,|y|^a)}
 = \int_D  |y|^a |u|^2 \, dx\, dy+\int_D |y|^a |\na u|^2 \, dx\, dy.
$$
\vspace{10pt}

We then introduce the energy functional:
$$\J(u) = \int_D \frac{1}{2} |y|^a |\na u|^2 \, dx\, dy$$
and the set
$$K_\eps =\left\{v\in W^{1,2}(D,|y|^a)\, ;\, v(x,0) \geq \vphi(x) \mbox{ for }x\in T_\eps(\omega)\, , \; v=g \mbox{ on } \Gamma \right\}.$$
It is readily seen that (\ref{eq:3}) is the Euler-Lagrange equation associated to the minimization problem: 
\begin{equation} \label{eq:inf}
\J (u^\eps) = \inf_{v \in K^\eps } \J (v) ,\quad \qquad u^\eps\in K_\eps.
\end{equation}
(Note that since $K_{\eps}$ is closed, convex and not empty, (\ref{eq:inf}) has a unique solution $u^{\eps}\in K_\eps$).
\vspace{10pt}

Finally, we notice (see \cite{CS}) that if $u(x,y)$ is the extension of  a function $f(x)$ as in (\ref{eq:extension}), then
$$ \int_{\RR^{n+1}_+} |y|^a |\na u|^2 \, dx\, dy = \int_{\RR^n} |\xi|^{2s} \big| \hat f(\xi)\big|^2\, d\xi = ||f||_{\stackrel{.}{H}^s(\RR^n)}.$$
In particular,  the minimization problem  (\ref{eq:inf}) is equivalent to the variational formulation of problem (\ref{eq:1}).

\vspace{20pt}

In this paper, we study the asymptotic behavior of the solutions of (\ref{eq:inf}) for any open subset $D$ of $\RR^{n+1}_+$.
The assumptions and the main result are made precise in the next section.
The proof of the main theorem, which is details in Section \ref{sec:thm}, relies on the construction of an appropriate corrector. This construction is detailed in Sections \ref{sec:balls} and \ref{sec:second}.
\vspace{10pt}

\section{Assumptions and Main result}
\subsection{The set $T_\eps$}
We consider  a probability space $(\Omega,\mathcal F, \mathcal P)$.
For all $\omega\in\Omega$,  the set $T_\eps(\omega)$ is given by:
$$T_\eps(\omega) = \bigcup_{k\in\ZZ^n} S_\eps(k,\omega)  ,$$
where the sets $S_\eps(k,\omega)\subset \RR^n$ satisfy the following assumptions:
\vspace{10pt}

\noindent {\bf Assumption 1:} 
{\it 
For all $k\in \ZZ^n$ and $\omega\in \Omega$, there exists $\gamma(k,\omega)$ (independent of $\eps$) such that
$$ \capa_s (S_\eps(k,\omega)) = \eps^n \, \gamma(k,\omega),$$
where $\capa_s (A)$ denotes the s-capacity of subset $A$ of $\RR^{n+1}$ (defined below). 

Moreover, we assume that
\begin{equation}\label{eq:M}
S_\eps(k,\omega) \subset B_{M\eps^\frac{n}{n-2s}}(\eps k) \mbox{ for all } k\in \ZZ^n\mbox{ a.e. } \omega \in \Omega,
\end{equation}
for some large constant $M$,
and that there exists a constant $\overline \gamma>0$ such that
\begin{equation}\label{eq:gb} 
 \gamma (k,\omega) \leq \overline \gamma \quad \mbox{ for all $k\in\ZZ^n$ and a.e. $\omega\in\Omega$}.
 \end{equation}
}

This first assumption defines the critical size of the set $T_\eps$. It will guarantee that $\capa_s(T_\eps) $ remains finite as $\eps$ goes to zero.
A natural definition for $s$-capacity of a subset $A$ of $\RR^n$ is the following:
$$ \capa_s (A) = \inf\left\{\int_{\RR^{n}}\!\!\! |\xi|^{2s} \big|\hat f(\xi)\big|^2\, d\xi\, ;\; f\in H_0^s(\RR^{n}),\; f(x)\geq 1 \mbox{ for }x\in A
\right\}.$$
Using the extension problem for the fractional Laplce operator (see \cite{CS} for details), an equivalent definition (up to a multiplicative constant) is given by
$$ \capa_s (A) = \inf\left\{\int_{\RR^{n+1}} \!\!\!\!\!\!\!\! y^a |\na h|^2\, dx\, dy\, ;\, h\!\in\! W_0^{1,2}(\RR^{n+1}_+\!\! ,|y|^a) ,\,   h(x,0)\geq 1 \mbox{, } x\in A
\right\}\! .$$
We will use this second definition in this paper.
If $B^n_r$ is a n-dimensional ball, then its $s$-capacity in $\RR^{n+1}$ is given by
$$ 
\capa_s (B^n_r)= c_{n+1-a} r^{n-1+a} = c_{n+2s} r^{n-2s}
$$ 
for some constant $c_k$.
Assumption {\bf 1} is thus  satisfied in particular if the sets 
$
 S_{\eps}(k,\omega) $
 are balls 
 centered on $\eps\ZZ^n$  with radius 
  $r(k,\omega) \eps^{\frac{n}{n-2s}} $.

\vspace{10pt}

{\bf Assumption 2:} 
{\it The process
$ \gamma : \ZZ^{n} \times \Omega \mapsto [0,\infty)$
is stationary ergodic:  There exists a family of measure-preserving transformations $\tau_k: \Omega\rightarrow \Omega$ satisfying
$$ \gamma(k+k',\omega)=\gamma(k,\tau_{k'}\omega)\quad \mbox{ for all } k,k'\in\ZZ^n\mbox{ and } \omega\in\Omega,$$
and such that 
if $A\subset \Omega$ and  $\tau_k A=A $ for all $k\in\ZZ^n$, then $P(A)=0$ or $P(A)=1$
(the only invariant set of positive measure is the whole set).
}

\vspace{15pt}

This second assumption is necessary to ensure that some averaging process occur as $\eps$ goes to zero (the hypothesis of stationarity is the most general extension of the notions of periodicity and almost periodicity for a function to have some self-averaging behavior). 
\vspace{10pt}

\vspace{15pt}

\subsection{Main result}
We are now ready to state our main result:
\begin{theorem}\label{thm:1}
Let $D$ be a open subset of $\RR^{n+1}_+$ ($n\geq 2$), denote 
$$\Sigma = D\cap\{y=0\},\qquad \Gamma = \pa D \cap \{y>0\}$$ 
and let $T_\eps(\omega)$ be a subset of $\Sigma $ satisfying Assumptions {\bf 1} and {\bf 2}  above.

There exists a constant $\alpha_0\geq 0$ such that
for any  $\vphi(x,y)\in \mathcal C^{1,1}(D)$ 
the solution $u^\eps(x,y,\omega)$ of
$$ \min \left\{ \frac{1}{2}\! \int_D \!\! y^a |\na v|^2 \, dx\, dy\, ;\,  v \in  \Wao, \, v(x,0)\geq\vphi(x,0) \mbox{ for }x\in T_\eps(\omega)\right\},
$$
converges $\Wa$-weak and almost surely $\omega\in\Omega$
to a function  $\ou(x,y)$  solution of
the following minimization problem
$$ \min \left\{\frac{1}{2}\! \int_D \!\! y^a |\na v|^2 \, dx\, dy +\frac{1}{2} \!\int_\Sigma \!\alpha_0 (v-\vphi)_-^2(x,0) \, dx \, ;\,  v \in  \Wao\right\},
$$
where $w_-=\max(0,-w).$

If, moreover,  there exists $\underline \gamma>0$ such that 
$\gamma (k,\omega) \geq \underline \gamma $  for all $k\in\ZZ^n$ and a.e. $\omega\in\Omega$,
then $\alpha_0>0$.
\end{theorem}
\vspace{10pt}

In particular the function $\ou(x,y)$ solves
$$\left\{
\begin{array}{ll}
 -\div(y^a \na \ou) = 0 & \mbox{ for }(x,y)\in D \\[5pt]
\ds \lim_{y\rightarrow 0} y^a \pa_y \ou (x,y)=  \alpha_0 (\ou-\vphi)_- (x,0) & \mbox{ for } x\in \Sigma \\[9pt]
\ou(x,y)=0& \mbox{ for } (x,y)\in \Gamma
\end{array}\right.
$$
 
\begin{remark} When $D$ is a bounded subset of $\RR^{n+1}_+$, the condition $ u \in  \Wao$ could easily be replaced by
$$u\in\Wa\, , \qquad u(x,y)=g(x,y)\mbox{ for }x\in\pa D\cap \{y>0\}$$
for some function $g(x,y)\in L^\infty(\pa D\cap \{y>0\})$.
\end{remark}

We stated Theorem \ref{thm:1} in its most general form. It contains the semipermeable membrane problem, as well as our original problem (\ref{eq:1}) with the fractional operator. More precisely, if we have $D=\RR^{n+1}_+$ and if we consider the trace  $v(x)=\overline u(x,0)$ in Theorem \ref{thm:1} we get:
\begin{corollary}
Let $T_\eps$ be a subset of $\RR^n$ ($n\geq 2$) satisfying Assumptions {\bf 1} and {\bf 2} above. 
There exists $\alpha_0\geq 0$ such that for any $\vphi(x)\in C^{1,1}(\RR^n)$,
the solution $v^\eps(x,\omega)$ of (\ref{eq:1}) converges, as $\eps\rightarrow 0$,  $H^s(D)$-weak and almost surely 
to a function $\overline v(x)$  solution of
$$
(-\Delta)^{s} v - \alpha_0 (v-\vphi)_- =0 .
$$
\end{corollary}

As in Cioranescu - Murat \cite{CM1,CM2} and Cafarelli-Mellet \cite{CM}, the proof of Theorem \ref{thm:1} relies on the construction of an appropriate corrector.
More precisely, we use the  following result:
\begin{proposition}\label{prop:1}
Let $T_\eps(\omega)$ be a subset of $\RR^n$ satisfying
 Assumptions {\bf 1} and~{\bf 2}  above. There exists  a non-negative constant $\alpha_0$ such that 
 for every bounded subset $D$ of $\RR^{n+1}_+$, 
there is a function $w^\eps_0(x,y,\omega)$ defined in $D$ and  satisfying 
\begin{eqnarray}
&&w^\eps (x,0) = 1\quad \mbox{ for  } x\in T_\eps(\omega)\cap (D\cap\{y=0\}) \\[8pt]
&&\|w^\eps\|_{L^\infty(D)}\leq C\label{eq:infty}
\\[8pt]
&&w^\eps  \longrightarrow 0\qquad \Wa \mbox{-weak}\mbox{ a.s. } \omega \in \Omega \label{eq:cv}
\end{eqnarray}
and
\begin{equation}\label{H1'}
\left\{
\begin{array}{l}
\mbox{For all sequences $v^\eps(x,y,\omega)$ satisfying:} \\[5pt]
\quad\qquad\left\{ \begin{array}{l}
v^\eps (x,0)\geq 0\quad \mbox{ for } x\in T_\eps(\omega)\cap\Sigma \\[4pt]
||v^\eps||_{L^\infty(D )} \leq C \\[4pt]
v^\eps\longrightarrow v\quad \mbox{ in } \Wa-\mbox{weak, } \mbox{a.s.}
\end{array}\right.\\[15pt]
\mbox{and for any $\phi\in\mathcal D(D)$ such that $\phi\geq 0$, we have:}\\[5pt]
\qquad\displaystyle\lim_{\eps\rightarrow 0} \int_D y^a \na w^\eps \cdot \na v^\eps \phi\, dx\, dy  \geq -\int_\Sigma \alpha_0 \, \phi \, dx \\[8pt]
\mbox{with equality if $v^\eps (x,0) = 0 \mbox{ for } x\in T_\eps\cap\Sigma$.}
\end{array}
\right.
\end{equation}
\end{proposition}

The proof of Proposition \ref{prop:1}  will occupy most of this paper.
We stress the fact that Assumptions {\bf 1} and {\bf 2} are sufficient but by no mean necessary to the proof of  this Proposition. Any set $T_\eps(\omega)$ such that Proposition \ref{prop:1} holds would be admissible for Theorem \ref{thm:1}.

The condition (\ref{H1'}) may seem rather obscure and the next Lemma will suggest a  nicer (but stronger) condition to replace it. 
However (\ref{H1'}) is the condition that appears naturally  in the proof of Theorem \ref{thm:1}.
\begin{lemma}
Let $D$ be a bounded subset of $\RR^{n+1}_+$, and 
assume that $w^\eps$ satisfies 
\begin{equation}\label{H1}
\left\{
\begin{array}{ll}
-\div(y^a\na w^\eps)  = 0  \quad & \mbox{ for } (x,y)\in D \\[5pt]
w^\eps (x,0) = 1\quad &\mbox{ for }x\in T_\eps(\omega)\cap\Sigma \\[5pt]
\lim_{y\rightarrow 0} y^a \pa_y w^\eps (x,y) = \alpha_0 \quad& \mbox{ for } x\in\Se\cap\Sigma\\[5pt]
\lim_{y\rightarrow 0} y^a \pa_y w^\eps (x,y) \leq 0 \quad& \mbox{ for }x\in T_\eps\cap\Sigma
 \end{array}
\right.
\end{equation}
together with (\ref{eq:infty}) and (\ref{eq:cv}).
Then (\ref{H1'}) holds.
\end{lemma}
This lemma also gives an indication of how to construct $w^\eps(x,y,\omega)$: We will look for a constant $\alpha_0$ such that the solution of (\ref{H1}) converges to zero in $\Wa$-weak.

{\bf Proof:}
Let $v^\eps\in L^\infty(D)\cap\Wa$ be such that $v^\eps(x,0)\geq 0$ on $T_\eps\cap\Sigma$ and let $\phi$ be a smooth test function with compact support in $D$.  Then, we have:
\begin{eqnarray*} 
0&=& \int_{D} \div(y^a \na w^\eps)  \phi\, v^\eps\, dx \, dy\\
&=&\!\!\!- \int_{D} y^ a \phi \na w^\eps\cdot \na v^\eps\, dx \, dy- \int_{D}y^a \na \phi \cdot \na w^\eps   v^\eps \, dx\, dy\\
&& \!\!\! - \int_{\Sigma\setminus T_\eps} \lim_{y\rightarrow 0} (y^a\pa_y w^\eps) \, v^\eps\, \phi\, dx- \int_{T_\eps} \lim_{y\rightarrow 0} (y^a\pa_y w^\eps) \, v^\eps\, \phi\, dx.
\end{eqnarray*}
Since $\lim_{y\rightarrow 0} y^a \pa_y w^\eps (x,y) \leq 0$ and $v^\eps(x,0)\geq 0$ for $x\in T_\eps$, we deduce:
\begin{eqnarray} 
\int_D y^a \phi  \na w^\eps \cdot\na v^\eps \, dx\,dy  & \geq & \!\!- \int_{D}y^a \na \phi \na w^\eps  v^\eps \, dx\, dy- \int_\Se \alpha_0 \, v^\eps\, \phi\, dx.\nonumber \\
 & \geq &\!\! - \int_{D}y^a \na \phi \na w^\eps  v^\eps \, dx\, dy- \int_\Sigma \alpha_0 \, v^\eps\, \phi\, dx. \label{eq:nnd}
\end{eqnarray}
with equality if $v^\eps(x,0)=0$ for $x\in T_\eps$.
In order to pass to the limit in (\ref{eq:nnd}), we note that we have the following convergences:
$$  w^\eps  \longrightarrow 0\qquad \Wa\mbox{-weak}\qquad\mbox{ a.s. } \omega \in \Omega,$$
and
$$  v^\eps  \longrightarrow v\qquad  L^2(D,|y|^a)\mbox{-strong} \qquad\mbox{ a.s. } \omega \in \Omega.$$
Hence the first term in the right hand side of (\ref{eq:nnd}) goes to zero. 
Moreover we have
$$  v^\eps (\cdot,0) \longrightarrow v(\cdot ,0) \qquad H^{s} (\Sigma)\mbox{-weak}  \mbox{ and }L^2(\Sigma)\mbox{-strong}\qquad\mbox{ a.s. } \omega \in \Omega,$$
so (\ref{eq:nnd}) gives
\begin{eqnarray*} 
\lim_{\eps\rightarrow 0}  \int_D y^a \phi  \na w^\eps \cdot\na v^\eps \, dx\,dy \geq - \int_\Sigma \alpha_0 \, v \, \phi\, dx.
\end{eqnarray*}
with equality if $v^\eps(x,0)=0$ for $x\in T_\eps$.
\qed

\subsection{Related problems}
Before turning to the proof of Theorem \ref{thm:1}, we briefly mention other results that follow from  Proposition \ref{prop:1}:
If we consider energy functionals of the form
$$\J(v) = \frac{1}{2}\int_D  y^a |\na u|^2 \, dx\, dy + \int_\Sigma u\,  h \, dx$$
for some  $h\in L^\infty(\Sigma)$, then a proof similar to that of Theorem \ref{thm:1} shows that
 the homogenization of the following equation
$$
\left\{
\begin{array}{ll}
v^\eps(x) \geq \vphi(x) & \mbox{ for } x\in T_\eps  \\[5pt]
(-\Delta)^{s} v^\eps \geq h(x) &  \mbox{ for } x\in \RR^n  \\[5pt]
(-\Delta)^{s} v^\eps = h(x) & \mbox{ for }x\in \RR^n\setminus T_\eps \mbox{ and on } T_\eps \mbox{ if } v^\eps  > \vphi
\end{array}
\right.
$$
leads to
$$
(-\Delta)^{s} v - \alpha_0 (v-\vphi)_- = h  \qquad  \mbox{ in } \RR^n  .
$$

More interestingly, we can replace the constrain $v^\eps \geq\vphi $ on $T_\eps$ by a Dirichlet condition of the form $v^\eps=0$
on $T_\eps$. 
This amounts to minimizing $\J(v)$ in the convex set
$$K_\eps = \{v\in W^{1,2}(\RR^{n+1}_+,|y|^a) \, ;\, v(x,0) =0 \mbox{ for }x\in T_\eps(\omega) \},$$
The corresponding Euler equation is 
$$
\left\{
\begin{array}{ll}
(-\Delta)^{s} v^\eps(x) = h(x) & \mbox{ for } x\in \RR^n\setminus T_\eps\\[4pt]
v^\eps(x)=0 &\mbox{ for } x\in   \pa T_\eps .
\end{array}
\right.
$$
We can then show that  the solution $v^\eps(x)$ converges to a function $v(x)$ solution of
$$
(-\Delta)^{s} v - \alpha_0 (v-\vphi) = h  \qquad  \mbox{ in } \RR^n .
$$

\vspace{10pt}

\section{Proof of Theorem \ref{thm:1}} \label{sec:thm}
In this section, we prove that Theorem \ref{thm:1} follows from  Proposition \ref{prop:1}.
For the sake of simplicity, we assume that $D$ is a bounded domain in $\RR^{n+1}_+$. This allows us to take the corrector $w^\eps(x,y,\omega)$ given by Proposition \ref{prop:1} and corresponding to the domain $D$. 
When $D$ is unbounded, we note that every integral involving $w^\eps$ is computed with a compactly supported test function $\phi$. We can thus use, for each of them, the corrector $w^\eps$ corresponding to the domain $\mbox{supp }\phi$. The final result is of course independent of $w^\eps$.
\vspace{10pt}

The maximum principle and the natural energy estimate easily give that
$u^\eps$ is bounded in $L^\infty(D)\cap\Wao$ almost surely. In particular, there exists a function $\ou(x,u,\omega)$ such that  
$$ u^\eps \longrightarrow \ou \qquad \Wao-\mbox{weak}\quad \mbox{a.e. } \omega\in\Omega.$$ 
In order to prove Theorem \ref{thm:1}, we have to show that
\begin{equation}\label{eq:Ja}
 \Ja (\ou) = \inf_{v\in \Wao} \Ja(v)\quad \mbox{a.e. } \omega\in\Omega
 \end{equation}
where $\Ja$ is
the energy associated to the limiting problem, given by:
$$ \Ja (v) = \frac{1}{2} \int_Dy^a |\na v|^2 \, dx\, dy +\frac{1}{2} \int_\Sigma \alpha_0 (u-\vphi)_-^2 \, dx.$$
Equality (\ref{eq:Ja}) will be a consequence of the following  lemmas:
\begin{lemma}\label{lem:1}
For any test function $\phi \in \mathcal D(D)$, we have
$$ \lim_{\eps\rightarrow 0} \int_D y^a |\na w^\eps |^2 \phi\, dx\, dy = \int _\Sigma \alpha_0 \phi\, dx.$$
\end{lemma}

\begin{lemma}\label{lem:2}
Let $u^\eps$ be a bounded sequence in $\Wa\cap L^\infty(D)$.
If 
$$u^\eps \rightharpoonup \ou \quad\mbox{ in  $\Wa$-weak,}$$ 
then
$$ \liminf_{\eps\rightarrow 0} \J(u^\eps) \geq \Ja(\ou).$$
\end{lemma}
\vspace{20pt}

\noindent{\bf Proof of Theorem \ref{thm:1}}:\\
For any $v\in \mathcal D(D)$, we consider the function $v+(v-\vphi)_- w^\eps$ (note that this function satisfies the obstacle constrain). Its energy is given by:
\begin{eqnarray*}
&& \J(v+(v-\vphi)_- w^\eps ) \\
&  &\quad =   \frac{1}{2} \int_D y^a \Big[ |\na v|^2 + |\na (v-\vphi)_-|^2 {w^\eps }^2+|(v-\vphi)_-|^2|\na w^\eps |^2\Big]\, dx\, dy \\
&&\qquad+\int_Dy^a \Big[ (v-\vphi)_-\na (v-\vphi)_- w^\eps \na w^\eps + \na v\na (v-\vphi)_- {w^\eps }\\
&& \qquad\qquad\qquad\qquad\qquad\qquad\qquad\qquad\qquad+\na v (v-\vphi)_-\na w^\eps \Big]\, dx\, dy.
\end{eqnarray*}
Lemma \ref{lem:1} and the weak convergence of $w^\eps $ to $0$ in $\Wa$ thus implies 
$$ \lim _{\eps\rightarrow 0}  \J (v+(v-\vphi)_- w^\eps ) = \Ja (v).$$

Morever, it is readily seen that the function $v+(v-\vphi)_- w^\eps $ belongs to $K_\eps$. 
Since $u^\eps$ minimizes $\J$ on $K_\eps$, we deduce
$$ \J(v+(v-\vphi)_- w^\eps ) \geq \J(u^\eps) ,
$$
and therefore
$$ \Ja(v) \geq \limsup_{\eps\rightarrow 0} \J(u^\eps)
\qquad\mbox{ for all $v\in\mathcal D(D)$.}$$
On the other hand, Lemma \ref{lem:2}  gives
$$ \liminf_{\eps\rightarrow 0} \J(u^\eps) \geq \Ja(\ou).$$
and so
$$ \Ja(\ou)\leq \Ja(v) \qquad \mbox{
for all } v\in \mathcal D(D).$$
Equality  (\ref{eq:Ja}) follows by a density argument.
\qed

\vspace{10pt}

{\bf Proof of Lemma \ref{lem:1}:} This first lemma is a straightforward consequence of (\ref{H1'}): If we take   $v^\eps = 1-w^\eps$, 
we  have  $v^\eps(x,0) =0$ for $x\in T_\eps$, $v^\eps(x,y)$ bounded in $L^\infty(D)$ and 
$  v^\eps(x,y)$ converges to $1$ in $\Wa \mbox{-weak}$, $L^2(D,|y|^a)\mbox{-strong,}$ and almost surely  $\omega \in \Omega.$
We can thus use (\ref{H1'}), which  implies
\begin{eqnarray*} 
- \int_D y^a \phi  \na w^\eps \cdot \na (1-w^\eps) \, dx\, dy   \longrightarrow  \int_\Sigma \alpha_0 \, \phi\, dx,
\end{eqnarray*}
and so
\begin{eqnarray*} 
 \int_{D} y^a \phi  |\na w^\eps |^2\, dx \, dy \longrightarrow  \int_\Sigma \alpha_0 \, \phi\, dx
\end{eqnarray*}
for all $\phi\in\mathcal D(D)$.
\qed

 \vspace{10pt}

{\bf Proof of Lemma \ref{lem:2}:} Following  Cioranescu-Murat (see \cite{CM2}, Proposition~3.1), we evaluate the quantity 
$$ \int_D |y|^a |\na (u^\eps - (z+(z-\vphi)_- w^\eps))|^2\, dx\, dy$$  
for some test function $z$ with compact support in $D$ and then take the limit as $\eps$ goes to zero.

Using (\ref{eq:cv}), we obtain:
\begin{eqnarray*}
\liminf_{\eps\rightarrow 0} \int_D y^a|\na u^\eps|^2\, dx\, dy & \geq &
2\int_D y^a \na\ou\cdot \na z\, dx\, dy -\int_D y^a |\na z|^2\, dx\, dy\\
& & + 2\lim_{\eps\rightarrow 0 } \int_Dy^a    (z-\vphi)_-\na u^\eps \cdot  \na w^\eps\, dx\, dy\\
& & - \lim_{\eps\rightarrow 0 }\int_D  y^a   (z-\vphi)_-^2 |\na w^\eps|^2\, dx\, dy.
\end{eqnarray*}
Lemma \ref{lem:1} yields
$$\lim_{\eps\rightarrow 0 }\int_D y^a    (z-\vphi)_-^2 |\na w^\eps|^2\, dx\, dy= \int_\Sigma \alpha_0  (z-\vphi)_-^2 \, dx.
$$
Property  (\ref{H1'}), together with the facts that  $u^\eps \in L^\infty(D)$ and $(u^\eps-\vphi)(x,0)\geq 0$ for $x\in T_\eps$, implies  
\begin{eqnarray*}
\lim_{\eps\rightarrow 0 } \int_D y^a(z-\vphi)_-\na u^\eps   \cdot  \na w^\eps 
& = & \lim_{\eps\rightarrow 0 } \int_D y^a (z-\vphi)_-\na (u^\eps-\vphi)  \cdot  \na w^\eps\, dx\, dy\\
&& +\lim_{\eps\rightarrow 0 } \int_Dy^a   (z-\vphi)_-\na \vphi  \cdot \na w^\eps\, dx\,dy\\
& \geq & - \int_\Sigma \alpha_0 (\overline u-\vphi)     (z-\vphi)_- \, dx .
\end{eqnarray*}
It follows  that for any test function $z\in\mathcal D(D)$ we have:
\begin{eqnarray*}
\liminf_{\eps\rightarrow 0} \int_D y^a|\na u^\eps|^2\, dx\, dy & \geq &
2\int_D y^a \na \ou\cdot \na z\, dx\, dy -\int_D y^a |\na z|^2\, dx\, dy\\
& & -2 \int_\Sigma \alpha_0 (\overline u-\vphi)     (z-\vphi)_- \, dx\\
& & -  \int_\Sigma \alpha_0  (z-\vphi)_-^2 \, dx.
\end{eqnarray*}
We can now take a sequence  $z_n$ that converges to $\overline u$ strongly in $\Wa$ and such that $z_n(\cdot ,0)$ converges to $\ou(\cdot,0)$ strongly in $L^2(\Sigma,|y|^a)$.
Using the fact that $ (\overline u-\vphi)     (\overline u-\vphi)_- = -(\overline u-\vphi)_-^2 $, we get
\begin{eqnarray*}
\liminf_{\eps\rightarrow 0} \int_D y^a|\na u^\eps|^2\, dx\, dy & \geq &
\int_D y^a |\na \ou |^2\, dx\, dy+  \int_\Sigma \alpha_0  (\ou-\vphi)_-^2 \, dx.
\end{eqnarray*}
which concludes the proof.\qed
\vspace{20pt}

\section{The auxiliary corrrector}\label{sec:balls}
\subsection{Notations and scheme of the proof}
We recall that
$$\RR^{n+1}_+ = \{(x,y)\in\RR^n\times\RR\, ;\, y\geq 0\},$$
and we fix a bounded domain $D\subset \RR^{n+1}_+$. 
For any $x_0\in\RR^n$ and $y_0>0$, we introduce the following notation for the Euclidian balls:
$$
\begin{array}{l}
B_r(x_0,y_0) = \left\{(x,y)\in \RR^{n+1}\, ;\,\left( |x-x_0|^2+|y-y_0|^2\right)^{1/2} \leq r \right\},\\[10pt]
B^{+}_r(x_0,0) =B_r(x_0,0)\cap \{y>0\}   ,  \\[8pt]
B^n_r(x_0) = \left\{x\in \RR^n\, ;\, |x-x_0|\leq r \right\}.
\end{array}
$$
\subsubsection{The fundamental solution}
We  recall  (see \cite{CSS} for details) that the function
$$ h(x,y) = \frac{\nu_{n+1+a}}{|x^2+y^2|^{\frac{n-1+a}{2}}}\qquad \mbox{with }\qquad \nu_k = \frac{\pi^{\frac k 2}\Gamma(\frac{k-1}{2})}{4},$$
solves
$$
\left\{
\begin{array}{l}
-\div(y^a\na  h)  (x,y) = 0 \quad \mbox{ for $y>0$ }\\[8pt]
\ds  \lim_{y\rightarrow 0} y^a\pa_y h(x,y) \longrightarrow -\delta(x),
\end{array}
\right.
$$
where $\delta(x)$ denotes the Dirac distribution centered at $0$ in $\RR^n$.
We also have
$$ \div (y^a \na h) = -\mu_{n,a}\delta(x,y)\qquad \mbox{in } \RR^{n+1}$$
where $\delta(x,y)$ denotes the Dirac distribution centered at $0$ in $\RR^{n+1}$ and for some  constant $\mu_{n,a}$.

\subsubsection{An auxiliary corrector}
One of the key point in the proof of Proposition \ref{prop:1} is to see that away from $\eps k$, the set $S_\eps(k,\omega)$ is equivalent to a (n+1)-dimensional ball.  More precisely,
we introduce the capacitary potential  $\vphi^\eps_{k}(x,y,\omega)$ associated to the set $S_\eps(k,\omega)$. 
It is defined by  the following minimization problem:
$$
\inf\left\{\int_{\RR^{n+1}} \!\!\!\! y^a |\na \vphi |^2\, dx\, dy\, ;\; \vphi\in W^{1,2}(\RR^{n+1}_+,|y|^a),\; \vphi(x,0)\geq 1 \;\forall x\in S_\eps(k,\omega)
\right\}.
$$
It is readily seen that, almost surely in $\omega$,   $\vphi^\eps_{k}(x,y,\omega)$ satisfies
\begin{equation}\label{eq:vphik}
\left\{ \!\! \begin{array}{ll} 
-\div(y^a \na \vphi^\eps_{k}) = 0   & \mbox{ for } (x,y) \in \RR^{n+1}_+\\[5pt]
\vphi^\eps_{k}(x,0) = 1\quad \qquad\qquad\qquad& \mbox{ for } x \in S_{\eps}(k,\omega) \\[5pt]
\lim_{y\rightarrow 0} y^a \pa_y \vphi^\eps_{k} (x,y) = 0& \mbox{ for } x \notin S_{\eps}(k,\omega) \\[5pt] 
\end{array}
\right.
\end{equation}
and by definition of the capacity as seen in the introduction, Assumption {\bf 1} yields
\begin{equation} \label{eq:gradg}
 \int_{\RR^{n+1}} y^a |\na \vphi^\eps_{k} |^2\, dx\, dy = \eps^n \gamma(k,\omega).
\end{equation}
Moreover, we have the following lemma (the proof of which is presented in Appendix \ref{app:capa}):
\begin{lemma}\label{lem:cap}
For any $\delta>0$, there exists $R_\delta$ such that
$$ \left| \vphi^\eps_{k}(x,y,\omega) - \eps^n \gamma (k,\omega) \frac{2}{\mu_{n,a}}h(x-\eps k,y)\right| \leq \delta \eps^n \gamma (k,\omega) \frac{2}{\mu_{n,a}} h(x-\eps k,y) \quad$$
for all $(x,y)$ such that $|(x-\eps k,y)| \geq \eps^{\frac{n}{n-1+a}} R_\delta$
and for all $\eps>0$.

Moreover, $R_\delta$ depends only on the constant $M$ appearing in Assumption~{\bf1} (in particular, $R_\delta$ is independent of $k$ and $\omega$).
\end{lemma}
This Lemma will play a fundamental role in the proof of Proposition \ref{prop:1} (see Section \ref{sec:second}).
It suggests that
at distance $\eps^{\frac{n}{n-1+a}}  R$ away from $\eps k$, the corrector $w^\eps$ should behave like the function 
$$ h_{k}^\eps (x,y,\omega) := \eps^n \gamma(k,\omega) \frac{2}{\mu_{n,a}} h(x-\eps k,y) .$$ 
For later use, we introduce the notation
$$a^\eps =\eps^{\frac{n}{n-1+a}} .$$
The first step in the proof, and the main goal of this section is to construct a function $\tw^\eps$ that would be a good approximation of $w^\eps$ away from  $\eps k$ and that behaves like $h_k^\eps$ at distance  $a^\eps R$ from $\eps k$
\vspace{10pt}

For that purpose, we introduce 
$$\widetilde D_\eps=D\setminus \bigcup_{k\in\ZZ^n}  B^+_{r(k,\omega) a^\eps }(\eps k), \quad \mbox{and  }\quad \widetilde \Sigma_\eps = \Sigma \setminus  B^n_{r(k,\omega) a^\eps }(\eps k),$$
where $r(k,\omega)$ is chosen in such a way that 
$h_{k}^\eps (x,y) =1$ on $\pa B^+_{r(k,\omega) a^\eps }(\eps k)$, i.e.
\begin{equation}\label{eq:r}
r(k,\omega) 
=
\left(\frac{2\nu_{n+1+a}}{\mu_{n,a}}\, \gamma(k,\omega)\right)^{1/(n-1+a)} .
\end{equation}
We will prove the following proposition:
\begin{proposition}\label{prop:aux}
There exist  a non-negative real number $\alpha_0$ (independent of the choice of $D$) and a function  $\tw^\eps(x,y,\omega)$ satisfying
\begin{equation}\label{H1bis}
\left\{
\begin{array}{ll}
\ds -\div(y^a\na  \tw^\eps)  = 0  \quad & \mbox{ for }(x,y)\in \widetilde D_\eps \\[8pt]
\ds \lim_{y\rightarrow 0} y^a\pa_y \tw^\eps(x,y) = \alpha_0 \quad & \mbox{ for  } x\in\widetilde \Sigma_\eps 
 \end{array}
\right.
\end{equation}
for almost all $\omega \in \Omega$, such that
\begin{equation}\label{eq:hkeps}
  \tw^\eps (x,y)\; = \; h_{k}^\eps (x,y) +o(1)\quad\mbox{ for } (x,y)\in  B^+_{\eps/2}(\eps k)\cap \widetilde D_\eps \mbox{ a.s. } \omega\in \Omega
\end{equation}
Moreover, we have:  
\item (i) $||\tw^\eps||_{L^\infty(\widetilde D_\eps) } \leq C$
\item (ii) $||\tw^\eps||_{L^2(\widetilde D_\eps) }\longrightarrow 0 $ as $\eps\rightarrow 0$.
\item (iii)$ ||\na \tw^\eps||_{L^2(\widetilde D_\eps) }\leq C$
\end{proposition}

The goal of this section is to establish Proposition \ref{prop:aux}. 
The main advantage of $ \tw^\eps$ over $w^\eps$ is that the former only depends on the capacity of $S_\eps(k,\omega)$. This explain why no assumptions are needed on the shape of $S_\eps(k,\omega)$.
In the last section of the paper  (Section \ref{sec:second}), we will see how to use both the functions $ \vphi^\eps_{k}$ (near $\eps k$) and  the corrector $\tw^\eps$ (at distance $a^\eps R$ of $\eps k$) in order to prove Proposition~\ref{prop:1}.
\vspace{10pt}

\subsubsection{Effective equation}
The main idea to prove  Proposition \ref{prop:aux} (and in particular (\ref{eq:hkeps})) makes use of the fact that $h_{k}^\eps (x,y,\omega)$
solves:
$$
\left\{
\begin{array}{ll}
-\div(y^a\na  h_k^\eps)  (x,y) = 0 \qquad &\mbox{ for $(x,y)\in\RR^{n+1}_+$ }\\[8pt]
\ds  \lim_{y\rightarrow 0} y^a\pa_y h_k^\eps(x,y) = -\eps^n\widetilde \gamma (k,\omega)  \delta(x-\eps k)& \mbox{ for $x\in\RR^{n}$ }
\end{array}
\right.
$$
with
$$\widetilde \gamma(k,\omega)=\gamma (k,\omega) \frac{2}{\mu_{n,a}}.$$
Proposition \ref{prop:aux} will thus be a consequence of the following proposition:
\begin{proposition}\label{lem:woe}
There exists $\alpha_0\geq 0$ such that
the solution $w^\eps_0(x,y,\omega)$ of
\begin{equation}\label{eq:wwo}
\left\{
\begin{array}{ll} 
-\div(y^a \na  w^\eps_0) = 0 &  \mbox{ for }(x,y)\in \RR^{n+1}_+ \\[8pt]
\ds \lim_{y\rightarrow 0} y^a\pa_y w^\eps_0 =  \alpha_0 - \!\!\! \sum_{k\in \ZZ^{n}\cap \Sigma} \eps^n\tgamma(k,\omega) \delta(x-\eps k)  \!\!\!& \mbox{ for }x\in \Sigma \\[8pt]
w^\eps_0(x,0) =0  & \mbox{ for } x\in\RR^n\setminus \Sigma
 \end{array}
 \right.
\end{equation}
satisfies:
\begin{equation}\label{eq:w0}
  w^\eps_0 (x,y)\; = \; h_{k}^\eps (x,y) +o(1)\quad\mbox{ for } (x,y)\in  B^+_{\eps/2}(\eps k)\cap D \mbox{ a.s. } \omega\in \Omega
\end{equation}
\end{proposition}

This proposition is the main step in the proof of Proposition \ref{prop:aux} and  its proof will occupy most of section.

\subsection{Proof of Proposition \ref{lem:woe}}
In order to prove Proposition \ref{lem:woe}, it is more convenient to work with the rescaled function 
\begin{equation}\label{rescale}
 v^\eps_0 (x,y,\omega) = \eps^{-1+a} w^\eps_0(\eps x,\eps y,\omega).
\end{equation}
Equation (\ref{eq:wwo}) then becomes:
\begin{equation}\label{vv0}
\left\{\begin{array}{ll}
- \div(y^a\na v^\eps_0) = 0\quad\qquad\quad & \mbox{ for } (x,y)\in \RR^{n+1}_+ \\[5pt] 
\ds \lim_{y\rightarrow 0} y^a \pa_y v^\eps_0(x,y) = \alpha_0 - \!\!\! \sum_{k\in\ZZ^{n}\cap D} \tgamma(k,\omega) \delta(x-k) \!\!\!\! &\mbox{ for } x\in \eps^{-1}\Sigma\\[8pt]
 v^\eps_0(x,0)=0 \qquad\qquad &\mbox{ for } x\in \RR^n\setminus \eps^{-1}\Sigma,
\end{array}
\right.
\end{equation}
and (\ref{eq:w0}) is equivalent to
$$  v^\eps_0 (x,y,\omega) \, =\, h_k(x,y,\omega) + o(\eps^{-1+a})\quad\mbox{ for } (x,y)\in  B^+_{1/2}
(\eps k)\cap \eps^{-1} D \mbox{ a.s. } \omega\in \Omega
$$ 
where
$$h_k(x,y):= \tgamma(k,\omega)\;  h(x-k,y) = \frac{ r(k,\omega)^{n-1+a}}{|(x-k)^2+y^2|^{(n-1+a)/2}}.
$$
Note that $h_k=\eps^{-1+a}$ on $\pa B_{\overline a^\eps r(k,\omega)}$ with $\overline a^\eps=\eps^\frac{1-a}{n-1+a}$.

In order to find the critical $\alpha_0$ for which the solution $v^\eps_0$ has the appropriate behavior near the lattice points $k\in\ZZ^n$, we follow the method developed by Caffarelli-Souganidis-Wang in \cite{CSW} and which was already the corner stone in \cite{CM}: We  introduce the following obstacle problem, for every open set $A\subset \RR^n$ and for every real number $\alpha\in\RR$:
\begin{equation}\label{vobstacle}
\left\{\begin{array}{ll}
v(x,0) \geq 0  &\mbox{ for } x\in \RR^{n} \\[5pt]
\ds \lim_{y\rightarrow \infty} v(x,y)=0&\mbox{ for } x\in \RR^{n} \\[5pt] 
- \div(y^a\na v^\eps) \geq 0\quad\qquad\quad& \mbox{ for } (x,y)\in \RR^{n+1}_+ \\[5pt] 
\ds \lim_{y\rightarrow 0} y^a \pa_y v(x,y) \leq \alpha - \sum_{k\in\ZZ^{n}\cap D} \tgamma(k,\omega) \delta(x-k)  &\mbox{ for } x\in A.
\end{array}
\right.
\end{equation}
We then define the smallest super-solution of the obstacle problem:
\begin{equation}\label{def:ov}
\ov_{\alpha,A}(x,y,\omega) = \inf\big\{v(x,y)\, ; \, v \mbox{ solution of  (\ref{vobstacle})} \big\}.
\end{equation}
It is readily seen that the function $\ov_{\alpha,A}$ satisfies
\begin{equation} \label{eq:v}
\left\{
\begin{array}{ll}
-\div(y^a \na \ov_{\alpha,A} ) = 0& \mbox{ for } (x,y)\in \RR^{n+1}_+  \\[8pt]
\ds\lim_{y\rightarrow 0 } y^a \pa_y \ov_{\alpha,A}(x,y) = \alpha\! -\!\!\!\!\! \sum_{k\in\ZZ^{n}\cap A}\!\! \tgamma(k,\omega) \delta(x-k) & \mbox{ for }x\in A\cap\{\ov_{\alpha,A}>0\}
\end{array}
\right.\end{equation}
and
\begin{equation}\label{eq:pay}
\lim_{y\rightarrow 0 } y^a \pa_y \ov_{\alpha,A}(x,y)  \geq0\qquad \mbox{ for }x\in A\cap\{\ov_{\alpha,A}=0\}.
\end{equation}

\vspace{20pt}

\begin{remark} The function 
\begin{eqnarray}\label{eq:hak}
\!\!\!\! h_{k,\alpha}(x,y) \!\!\!\! & = &\!\!\!\!  h_k(x,y)  - \alpha \int_{B^n_1(k)} \frac{\nu_{n+1+a}}{(|x-x'|^2+y^2)^\frac{n-1+a}{2}} dx' 
\end{eqnarray}
satisfies
$$\left\{
\begin{array}{ll}
-\div(y^a\na  h_{k,\alpha} )= 0&  \mbox{ for  }x\in \RR^{n+1}_+   \\[8pt]
\ds \lim_{y\rightarrow 0} y^a\pa_yh_{k,\alpha}(x,y) \longrightarrow \alpha  -\tgamma(k,\omega) \delta(x-k)&  \mbox{ for } x\in B^n_1(k).
\end{array}
\right.$$
It is radially symmetric around $k$ and 
$\sup_{|x|=1, \; y>0} h_{\alpha,k}(x,y)\leq r^{n-1+a}.$
In particular, the maximum principle and  (\ref{eq:v}) implies that
if $B^n_1(k)\subset A$, then:
\begin{equation}\label{eq:ineqw} 
\ov_{\alpha,A}(x,y,\omega)  \geq 
\ds h_{\alpha,k}(x,y,\omega) -r^{n-1+a} \qquad \mbox{ for } (x,y)\in  B^+_1(k), \mbox{ a.s.}
\end{equation}
\end{remark}

\vspace{10pt}

We now want to show that there exists a critical $\alpha_0$ such that the followings hold:
\begin{enumerate}
\item  The solution of the obstacle problem $\ov_{\alpha,A}(x,y,\omega)$ behaves like $h_{\alpha,k}(x,y,\omega)$ near any point $k\in A \cap\ZZ^n$.
\item The solution of (\ref{vv0}) is not  far from $\ov_{\alpha,A}$.
\end{enumerate}

For that purpose, we introduce the following quantity, which measures the size of the contact set along the boundary $\{y=0\}$: 
$$
\om_{\alpha}(A,\omega) = |\{x\in A\, ;\, \ov _{\alpha,A}(x,0,\omega)= 0 \}| 
$$
where $|A|$ denotes the Lebesgue measure of a set $A$ in $\RR^n$.
\vspace{10pt}

The starting point of the proof is the following lemma:
\begin{lemma}\label{lem:ergodic}
The random variable $\om_{\alpha}$ is subadditive, 
and the process 
$$T_k m(A,\omega) = m(k+A,\omega)$$
has the same distribution for all $k\in \ZZ^n$.
\end{lemma}
{\bf Proof of  Lemma \ref{lem:ergodic}:}
Assume that the finite family  of sets $(A_i)_{i\in I}$ is such that 
$$ 
\begin{array}{l}
A_i \subset A \qquad \mbox{ for all } i\in I \\
A_i\cap A_j = \emptyset \quad  \mbox{ for all } i\neq j \\
|A-\cup_{i\in I} A_i| = 0
\end{array}
$$ 
then $\ov_{\alpha,A}$ is admissible for each $A_i$, and so $\ov_{\alpha,A_i} \leq u_{\alpha,A}$.
It follows that 
$$ \{\ov_{\alpha,A}(\cdot,0,\omega)=0\}\cap A_i \subset\{\ov_{\alpha,A_i}(\cdot,0,\omega)=0\}$$
and so
\begin{eqnarray*}
\om_{\alpha}(A,\omega) & = &\sum_{i\in I} |\{\ov_{\alpha,A}(\cdot,0,\omega)=0\}\cap A_i | \\
&\leq&  \sum_{i\in I}|\{\ov_{\alpha,A_i}(\cdot,0,\omega)=0\}| =   \sum_{i\in I}\om_{\alpha}(A_i,\omega),
\end{eqnarray*}
which gives the subadditive property.
Assumption {\bf 2} then yields
$$ T_k m(A,\omega) = m(A,\tau_k \omega)$$
which gives the last assertion of the lemma.
\qed
\vspace{20pt}

Since 
$ \om_{\alpha}(A,\omega) \leq |A|,$ and thanks to the ergodicity of the transformations $\tau_k$,
it follows from the subadditive ergodic theorem (see \cite{DM}) that for each $\alpha$, there exists a constant $\oell(\alpha)$ such that
$$ \lim_{t\rightarrow \infty} \frac{\om_{\alpha}( B_t(0),\omega)}{|B_t(0)|} = \oell (\alpha)\quad  \mbox{ a.s., } \quad $$
where $B_t(0)$ denotes the ball centered at the origin with radius $t$.
Note that the limit exists and is the same if instead of $B_t(0)$, we use cubes or balls  centered at $tx_0$ for some $x_0$.

If we scale back and consider the function 
$$ \overline  w^\eps_\alpha (x,y,\omega)= \eps^{1-a} \; \ov_{\alpha,B_{\eps^{-1}}(\eps^{-1}x_0)}(x/\eps,y/\eps,\omega) ,\qquad \mbox{ in } B_1(x_0),$$
we deduce
$$ \lim_{\eps\rightarrow 0} \frac{|\{ x\, ;\, \overline w^\eps_{\alpha}(x,0,\omega)=0\}|}{|B_1|} = \oell (\alpha)\quad  \mbox{ a.s. }
$$

\vspace{10pt}

The next lemma summarizes the properties of  $\oell(\alpha)$:
\begin{lemma}\label{lem:alpha}
\item[(i)] $\oell(\alpha)$ is a nondecreasing functions of $\alpha$.
\item[(ii)] If $\alpha<0$, then $\oell (\alpha) = 0$. Moreover, if the  $\gamma(k,\omega)$ are bounded from below, then $\oell (\alpha) = 0$ for  $\alpha $  positive small enough ($0<\alpha < C(\underline \gamma)$).
\item[(iii)] If $\alpha$ is large enough ($\alpha \geq  C(\overline \gamma)$), then 
 $\oell (\alpha) > 0$.
\end{lemma}

The proof of this Lemma is rather technical and of little interest. It is presented in full details in  Appendix \ref{sec:alpha}.
Using Lemma \ref{lem:alpha},  we can define
$$\alpha_{0} = \sup \{\alpha\, ;\, \oell(\alpha) = 0 \}.$$
We observe that $\alpha_0$ is finite  (Lemma \ref{lem:alpha} (iii)) and that $\alpha_0$ is non negative (Lemma \ref{lem:alpha} (ii)). Moreover, $\alpha_0$ is  strictly positive if the $\gamma(k,\omega)$ are bounded from below almost surely by a positive constant.

\vspace{10pt}

We now fix a bounded subset $A$ of $\RR^n$ and 
we  denote by 
\begin{equation}\label{eq:ovae}
\ovae(x,y,\omega) =\ov_{\alpha,\eps^{-1}A}(x,y,\omega)
\end{equation}
the solutions of (\ref{def:ov}) corresponding to $\eps^{-1}A$. 
We also introduce  the rescaled function
$$ \wae (x,y,\omega)= \eps^{1-a} \; \ovae (x/\eps,y/\eps,\omega) .$$
\vspace{20pt}

In order to complete the proof of Proposition \ref{lem:woe}, we are first going to prove that $ \wae$ satisfies inequality (\ref{eq:w0}),
and then that the solution $w^\eps_0$ of (\ref{eq:wwo}) behaves like  $ \wae$.

We recall the definition of $h_{\alpha,k}$:
$$h_{\alpha,k} (x,y) = \frac{r(k)^{n-1+a} }{(|x-k|^2+y^2)^\frac{n-1+a}{2}} -\alpha \int_{B^n_1(k)} \frac{\nu_{n+1+a}}{(|x-x'|^2+y^2)^{\frac{n-1+a}{2}}}\, dx',
$$
and we introduce the scaled function
$$ 
h_{\alpha,k}^\eps (x,y):=\eps^{1-a} h_{\alpha,k}(x/\eps,y/\eps) .
$$
Note that when  $(x,y)\in \pa B^+_{ \overline a^\eps r(k,\omega)}(k)$, then
$$ h_{\alpha,k}(x,y)= 
  \eps^{-1+a} - \alpha \int_{B^n_1(0)} \frac{\nu_{n+1+a}}{(|x-x'|^2+y^2)^{\frac{n-1+a}{2}}}\, dx'
$$
(we recall that $\overline a^\eps=  \eps^{\frac{1-a}{n-1+a}} $).

We then have the following lemma:
\begin{lemma}\label{lem:1.2}
\item[(i)] For every $\alpha $ and for every $k\in\ZZ^n\cap A$, we have 
$$ \ovae(x,y)\geq
\ds h_{\alpha,k}(x,y) -r^{n-1+a} \quad \mbox{ for }(x,y)\in  B^+_{1}(k) \mbox{ a. s. }
$$
\item[(ii)] For every $\alpha > \alpha_0$ and  every $k\in\ZZ^n\cap A$, we have 
$$ \ovae(x,y)\leq h_{\alpha,k}(x,y)+ o(\eps^{-1+a})\quad \mbox{ for }(x,y)\in  B^+_{1/2}(k) \mbox{ a. s. }
 $$ 
\end{lemma}
We deduce:
\begin{corollary}
\item[(i)] For every $\alpha $ and every $k\in\ZZ^n\cap A$ such that $r(k,\omega)>0$, we have 
$$ \ovae(x,y)\geq \eps^{-1+a}+o(1) \quad \mbox{ for } (x,y)\in \pa B^+_{\overline r(k,\omega) \overline a^\eps }(k) \quad\mbox{ a.e. } \omega\in\Omega$$
and so
$$ \wae (x,y )\geq 1+ o(\eps^{1-a}) \quad \mbox{ for } (x,y)\in \pa B^+_{ r(k,\omega) a^\eps  }(k) \quad\mbox{ a.e. } \omega\in\Omega$$
for all $\alpha$. 
\item[(ii)] For every $\alpha > \alpha_0$ and every $k\in\ZZ^n\cap A$, we have 
$$ \ovae(x,y)\leq \eps^{-1+a}+o(\eps^{-1+a})\quad \mbox{ for }(x,y)\in \pa B^+_{r(k,\omega)\overline a^\eps}(k) \quad\mbox{ a.e. } \omega\in\Omega$$
and so
$$ \wae(x,y) \leq 1+ o(1)\mbox{ for }(x,y)\in \pa B^+_{ r(k,\omega)a^\eps}(k) \quad\mbox{ a.e. } \omega\in\Omega$$
\end{corollary}

\vspace{10pt}

\noindent {\bf Proof of Lemma \ref{lem:1.2}:}\\
\noindent (i) This is an immediate consequence of (\ref{eq:ineqw}).
\vspace{10pt}

\noindent (ii) The proof of (ii) is more delicate and is split in several steps.\\
{\bf Preliminary:}
First of all since $A$ is bounded, we have
$ A \subset B^n_R(x_0)$ for some $R$.
Without loss of generality, we can  always assume that $B^n_R(x_0)=B^n_1(0)$.
If we consider $$
\vae(x,y,\omega) =\ov_{\alpha,\eps^{-1}B^n_1}(x,y,\omega),
$$
the solution of (\ref{def:ov}) corresponding to  $A=B^n_{\eps^{-1}}(0)$,
it is readily seen that 
$$\ovae(x,y,\omega) \leq \vae(x,y,\omega)  \qquad \mbox{for all } (x,y) \in \RR^{n+1}_+\, \; \mbox{ a.e. }\omega\in\Omega.$$
It is thus enough to prove (ii) for $\vae$.
\vspace{10pt}

In the sequel, we will need the following consequence of Lemma \ref{lem:ergodic} (see \cite{CSW} for the proof):
\begin{lemma}\label{lem:local}
For any  ball $B^n_r(x_0)\in B^n_1(0)$, the following limit holds, a.s. in $\omega$:
$$\lim_{\eps\rightarrow 0} \frac{| \{\vae (x,0,\omega)= 0 \}\cap B^n_{\eps^{-1}r}(\eps^{-1}x_1)  | }{|B^n_{\eps^{-1}r}|} = \oell(\alpha)
$$
\end{lemma}
\vspace{10pt}

\noindent {\bf Step 1:} We now start the proof:
For any $\delta>0$, we can cover $B^n_{\eps^{-1}}$ by a finite number $N$ ($\leq C\delta^{-n}$) of balls $B^n_i=B^n_{\delta\eps^{-1}}(\eps^{-1}x_i)$ with radius $\delta \eps^{-1}$ and center $\eps^{-1}x_i$.
Since $\alpha >\alpha_0$, we have $\oell(\alpha)>0$. By Lemma \ref{lem:local}, we deduce that for every $i$,  there exists $\eps_i$ such that if $\eps\leq\eps_i$, then
$$| \{\vae (x,0,\omega)= 0 \}\cap B^n_{i}  | >0\quad \mbox{ a.s. } \omega.$$
In particular, if $\eps \leq \inf \eps_i$, then $\vae(x'_i,0)=0$ for some $x'_i$ in $B^n_i$ a.s. $\omega\in\Omega$.

Introducing $B_i=B_{\delta\eps^{-1}}(\eps^{-1}x_i)$ the $n+1$ dimensional ball with same radius and same center as $B^n_i$,
we now have to show that  $\vae$ remains small in each $B^+_i$ as long as we stay away from the lattice points $k\in\ZZ^n$.
More precisely, we want to show that
$$ \sup_{ \cup_{k\in\ZZ^n} B^+_1(k)\setminus B^+_{1/4}(k)} \vae(x,y) \leq C\delta^{1-a} \eps^{-1+a}.$$

\vspace{20pt}

\noindent {\bf Step 2:} 
Let $\eta(x)$ be a nonnegative function defined in $\RR^n$ such that $0\leq \eta(x)\leq 1$ for all $x$,  $\eta(x)=1$ in $ B_{1/8}$ and $\eta = 0$ in $\RR^{n} \setminus B_{1/4}$.
We then consider the function
$$u=\vae \star_x \eta$$
where $\star_x$ indicates the convolution in $\RR^n$ with respect to the $x$-variable. 
The function $u(x,y)$ is nonnegative on $2 B_i^+$ and satisfies 
\begin{equation}\label{eq:lap}
\left\{
\begin{array}{ll}
\div(y^a\na u) = 0 & \mbox{ for }(x,y)\in 2 B^+_i\\[3pt]
 -C\leq \lim_{y\rightarrow 0} y^a \pa_y  u(x,y) \leq C &\mbox{ for }x\in  2 B^n_i
\end{array}
\right.
 \end{equation}
where $C$ is a universal constant depending only on $n$, $\overline r$ and $\alpha$.
We deduce:
\begin{lemma}\label{lem:harnack}
There exists a universal constant $C$ such that
$$\sup_{ B_i} u \leq C \inf_{ B_i}   u + C \delta^{1-a} \eps^{-1+a} .$$
\end{lemma}
{\bf Proof:} We write $u=u_1+u_2$ where $u_1$ and $u_2$ are two functions solution of $\div(y^a \na u_i )= 0$  in  $2  B^+_i$ and satisfying
$$
\left\{
\begin{array}{ll} 
\ds \lim_{y\rightarrow 0}\,  y^a \pa_y u_1(x,y) =  \lim_{y\rightarrow 0}\, y^a \pa_y  u (x,y) & \mbox{ for }  x\in 2B^n_i,\qquad \\[5pt]
u_1(x,y) = 0 & \mbox{ for } (x,y)\in\pa (2  B^+_i)\cap \{y>0\}
\end{array}
\right.
$$
and 
$$
\left\{
\begin{array}{ll}
\ds  \lim_{y\rightarrow 0} y^a \pa_y u_2(x,y) = 0&  \mbox{ for }x\in 2B^n_i,\qquad \\[5pt]
u_2 (x,y)= u(x,y) & \mbox{ for }(x,y)\in \pa( 2  B^+_i)\cap \{y>0\}.
\end{array}
\right.
$$
The maximum principle and the fact that $B_i$ has radius $\delta \eps^{-1}$ yield:
\begin{eqnarray*} 
|u_1(x,y)| &\leq & C ( (2\delta \eps^{-1})^{1-a}-  y^{1-a}) \\
&\leq& C ( \delta \eps^{-1})^{1-a}
\end{eqnarray*}
for all $(x,y)\in 2 B^+_i.$
On the other hand, boundary Harnack inequality  for degenerate elliptic equation (see \cite{FKS}) implies
$$\sup_{ B_i} u_2 \leq C \inf_{ B_i}   u_2 .$$
The Lemma follows easily.\qed

\vspace{20pt}

For the next step, we will need the following lemma:
\begin{lemma}\label{lem:average}
If $v$ satisfies 
$$\div( y^a \na v)=0 \quad \mbox{ in }  B^+_r(x_0,0) $$
and
 $$\lim_{y\rightarrow 0}y^a\pa_y v(x,y) \leq \alpha\quad \mbox{ for }x\in B^n_r(x_0),$$ 
 then
$$ \frac{2}{\omega_{n+a} r^{n+a}}\int_{B^+_r(x_0,0)} |y|^a  v (x,y)\, dx\, dy \leq v (x_0,0) + \alpha C(n) r^{1-a}$$
where $C(n)$ is a universal constant and 
 $\omega_{n+a} = \int_{ B_1(x_0,0)} |y|^a\, dx\, dy$.
\end{lemma}
{\bf Proof:} The function  
$w(x,y) = v(x,y)+\alpha \int_{B^n_r(x_0)} \frac{C_{n+1+a}}{(|x-x'|^2+y^2)^{\frac{n-1+a}{2}}}\, dx'$ 
satisfies
$$ \div(y^a \na w) =0 \quad\mbox{ and } \quad \lim_{y\rightarrow 0} y^a \pa_y w \leq 0.$$
Proceeding as in \cite{CS}, we now reflect $w$ about the plane $\{y=0\}$. The function 
$$ \overline w(x,y) = \left\{
\begin{array}{ll}
w(x,y) &\mbox{ if } y>0\\[5pt] 
w(x,-y) &\mbox{ if } y<0 
\end{array}
\right.
$$ 
is now defined in the whole space $\RR^{n+1}$ and it satisfies 
$$\div ( |y|^a \na \overline w)\leq 0 \quad \mbox{ in }  B_r(x_0,0).$$
We can thus use the mean value formula (see \cite{CSS}):
\begin{eqnarray*}
&& \frac{1}{\omega_{n+a} r^{n+a}}\int_{B_r(x_0,0)} |y|^a\overline  w (x,y)\, dx\, dy \\
&& \qquad \qquad \qquad \qquad \leq \overline w(x_0,0)\\
&& \qquad \qquad \qquad \qquad \leq w(x_0,0)\\
&& \qquad \qquad \qquad \qquad \leq  v(x_0,0)+\alpha \int_{B^n_r(x_0)} \frac{C_{n+1+a}}{|x_0-x'|^{n-1+a}}\, dx'.
\end{eqnarray*}
Since $\alpha\geq 0$, we see that $v\leq w$ and so
\begin{eqnarray*}
\frac{2}{\omega_{n+a} r^{n+a}}\int_{ B_r^+(x_0,0)}\!\!\!\! y^a  v (x,y)\, dx\, dy & \leq & 
\frac{1}{\omega_{n+a} r^{n+a}}\int_{ B_r (x_0,0)} \!\!\!\! |y|^a \overline w (x,y)\, dx\, dy \\
\end{eqnarray*}
Moreover, we have 
$$\int_{B^n_r(x_0)} \frac{C_{n+1+a}}{|x_0-x'|^{n-1+a}}\, dx' = \int_{B^n_r(0)} \frac{C_{n+1+a}}{|z|^{n-1+a}}\, dz = C(n+a) r^{1-a}, $$ 
hence the lemma.
\qed
\vspace{20pt}

\noindent {\bf Step 3:}
We have $\vae(x'_i,0)=0$ and $\lim_{y\rightarrow 0}y^a \pa_y\vae(x,y)  \leq \alpha$ for  $x\in B_{1/2}(x'_i)$.
Lemma \ref{lem:average} thus applies and yields:
\begin{equation}\label{eq:mf1} 
\int_{ B^+_{1/2}(x'_i,0)}\!\!\!\!\! |y|^a \vae (x,y)\, dx\, dy \leq C (\vae (x'_i,0) + \alpha )  \leq  C (\alpha,n+a).
\end{equation}
We want to deduce an upper bound on  $u$ in $B_i$. 
Since $u\geq 0$, we note that
$$\int_0^{1/4}\tau^a u(x,\tau) \, d\tau \geq \left(\inf_{\tau\in[0,1/4]} u \right) \int_0^{1/4} \tau^a\, d\tau.$$
Then, using the definition of $u$ (and the fact that  $\eta(x)=0 $ outside $B^n_{1/4}(x)$), we deduce:
\begin{eqnarray*}
\inf_{B^+_{1/4}(x'_i,0)} u & \leq &C \inf_{x} \int_0^{1/4}\tau^a u(x,\tau) \, d\tau\\
 & \leq & C \inf_{x}\int_0^{1/4} \int_{B^n_{1/4}(x)} \tau^a  \vae(\xi,\tau) \, d\xi\, d\tau\\
& \leq & C \int_{B_{1/2}(x_i',0)} \tau^a  \vae(\xi,\tau) \, d\xi\, d\tau,
\end{eqnarray*}
Which, together with (\ref{eq:mf1}) yields:
\begin{equation}\label{eq:infu} 
\inf_{ B^+_{1/4}(x'_i,0)} u \leq  C (\alpha,n).
\end{equation}
Using Lemma \ref{lem:harnack} we see that for every $\delta$ and for $\eps$ small enough, we have:
\begin{equation}\label{eq:step3}
\sup_{ B_i} u \leq C \inf_{ B_i}   u + C \delta^{1-a} \eps^{-1+a} \leq C(\alpha,n)+C \delta^{1-a} \eps^{-1+a}\leq C \delta^{1-a} \eps^{-1+a}.
\end{equation}

\vspace{20pt}

\noindent {\bf Step 4:}
We now want to use (\ref{eq:step3}) to get an upper bound on $\vae$. For that purpose, we note that 
$\lim_{y\rightarrow 0}y^a \pa_y \vae\geq 0$ in $B_i\setminus \cap_{k\in\ZZ^n}\{k\}$, and so a proof similar to that of Lemma \ref{lem:average} yields
\begin{equation}\label{eq:mf2}
\vae (x,y)\leq C_{n+a}\int_{B_{1/8}^+(x,y)}|\tau|^a \vae (\xi,\tau)\, d\xi \, d\tau  
\end{equation}
for all $(x,y)\in  B_i\setminus \cap_{k\in\ZZ^n} \overline B_{1/4}(k)$.

Inequality (\ref{eq:mf2}) and the definition of $u(x,y)$ yield that for all $(x,y)$ in $B_i\setminus \cap_{k\in\ZZ^n}  B_{1/4}(k)$, we have:
\begin{eqnarray*}
\vae (x,y)& \leq &   C_{n+a}\int_{y-1/8}^{y+1/8} \int_{ B^n_{1/8}(x)}|\tau|^a \vae (\xi,\tau)\, d\xi \, d\tau  \\
& \leq &     C_{n+a}\int_{y-1/8}^{y+1/8}  |\tau|^a u  (x,\tau)\,  \, d\tau  \\
& \leq &  C(n+a)|y|^{1+a} \sup_{ B_i} u  .
\end{eqnarray*}
Inequality  (\ref{eq:step3}) therefore implies
 \begin{equation}\label{eq:step4}
 \sup_{(x,y)\in  \cup_{k\in\ZZ^n} B^+_1(k)\setminus B^+_{1/4}(k)} \vae(x,y) \leq  C \delta^{1-a} \eps^{-1+a}.
 \end{equation}
\vspace{20pt}

 {\bf Step 5:}
In order to complete the proof of the lemma, we only have to notice that since $\inf_{\pa   B_{1/2}} h_{\alpha,k}(x,y) \geq -C \alpha$, (\ref{eq:step4}) and 
the definition of $\vae$ imply
$$
\vae(x,y) \leq h_{\alpha,k}(x,y) +  C \delta^{1-a} \eps^{-1+a}\qquad \mbox{ in }  B_{1/2}(k)
$$
for all $k\in\ZZ^n$. 
\qed
\vspace{30pt}

This conclude the proof of Lemma \ref{lem:1.2}, and 
we are now in position to complete the proof of Propositions \ref{lem:woe}.

\vspace{10pt}

\noindent {\bf Proof of Proposition \ref{lem:woe}.} \\
For every $\alpha$, we denote by $\overline v^\eps_\alpha$  the solution of the obstacle problem (\ref{vobstacle}) corresponding to $A=\eps^{-1}\Sigma$:
$$ \overline v^\eps_\alpha (x,y,\omega)=  \; \ov_{\alpha,\eps^{-1}\Sigma}(x,y,\omega) ,$$
and by
$\overline w^\eps_\alpha$  the rescaled function:
$$ \overline w^\eps_\alpha (x,y,\omega)= \eps^{1-a} \; \ov_{\alpha,\eps^{-1}\Sigma}(x/\eps,y/\eps,\omega) .$$
We recall that $w^\eps_0$ is solution of
$$
\left\{
\begin{array}{ll} 
-\div(y^a\na  w^\eps_0) = 0 &  \mbox{ for }(x,y)\in \RR^{n+1}_+ \\[3pt]
\ds \lim_{y\rightarrow 0} y^a \pa_y w^\eps_0(x,y) =  \alpha_0 - \sum_{k\in \ZZ^{n}\cap D} \tgamma(k,\omega) \delta(x-\eps k) & \mbox{ for } x\in\Sigma \\[3pt]
w^\eps_0(x,0) =0  & \mbox{ for } x\in\RR^n\setminus \Sigma
 \end{array}
 \right.
$$
In order to prove Proposition \ref{lem:woe}, we have to establish (\ref{eq:w0}). This is done in two steps using the properties of the function $\overline w^\eps_\alpha$:

\begin{enumerate}
\item For every $\alpha> \alpha_0$, we have
$ \div(y^a\na (w^\eps_0-w^\eps_\alpha)) =0 $ for $(x,y)\in\RR^{n+1}_+$,
$\ds\lim_{y\rightarrow 0}y^a \pa_y (w^\eps_0-w^\eps_\alpha) \geq \alpha_0-\alpha $ on $\Sigma$ and $(w^\eps_0-w^\eps_\alpha)(x,0)=0$ on $\RR^n\setminus \Sigma$. We deduce
$$ w^\eps_0(x_0,y_0)-w^\eps_\alpha(x_0,y_0) \leq \int_{\Sigma}  \frac{ \alpha_0-\alpha}{|(x_0-x)^2+y_0^2|^\frac{n-1+a}{2}}\, dx,$$
and therefore
$$ \sup_ {(x,y)\in\RR^{N+1}_+} \left( w^\eps_0(x,y)-w^\eps_\alpha (x,y) \right) \leq 
C|\Sigma  |^{\frac{1-a}{n+1}} \rho_\Sigma^{\frac{1-a}{n+1}}   \,|\alpha-\alpha_0| 
$$
with
$$
\rho_\Sigma  = \inf \{\rho\, ;\, \Sigma \subset B_\rho\}.
$$
In particular, we thus have
$$ w^\eps_0(x,y)\leq w^\eps_\alpha(x,y) +O(\alpha-\alpha_0) \qquad \mbox{ for } (x,y)\in\RR^{n+1}_+,$$
and Lemma \ref{lem:1.2} (ii) (since $\alpha> \alpha_0$) yields:
$$ w^\eps_0(x,y) \leq h_{\alpha,k}^\eps (x,y) +O(\alpha-\alpha_0)+o(1)\quad  \mbox{ for } (x,y)\in  B_{\eps/2}(\eps k) \mbox{ a.s. } $$ 
(Note that this argument  shows the continuity of $w^\eps_\alpha$ with respect to $\alpha$).

\item Similarly, we observe that for $\alpha \leq \alpha_0$, we have $\div(y^a\na (w^\eps_\alpha-w^\eps_0)) =0$ for $(x,y)\in\RR^{n+1}_+$, $(w^\eps_\alpha-w^\eps_0)(x,0)=0$ for $x\in\RR^n\setminus \Sigma$ and 
$$\lim_{y\rightarrow 0} y^a\pa_y  (w^\eps_\alpha-w^\eps_0)(x,y) \geq \alpha-\alpha_0 -\alpha 1_{\{w^\eps_\alpha = 0\}\cap \Sigma} \qquad
\mbox{ for $x\in\Sigma$}.$$

Proceeding as before, we deduce:
\begin{eqnarray*}
 \sup_{\RR^{N+1}_+}  \left( w^\eps_\alpha-w^\eps_0  \right)&  \leq &  C\rho_\Sigma^{\frac{1-a}{n+1}} \left[|\Sigma |^{\frac{1-a}{n+1}}(\alpha_0-\alpha)\right.\\
 && \qquad\qquad\qquad \left.+ C \alpha |\{w^\eps_\alpha (x,0)= 0\}\cap\Sigma|^{\frac{1-a}{n+1}}\right].
 \end{eqnarray*}
So Lemma \ref{lem:1.2} (i) yields
$$ w^\eps_0(x,y) \geq h_{\alpha,k}^\eps(x,y) - o(\eps) - O(\alpha_0-\alpha) - C \alpha |\{w^\eps_\alpha(x,0) = 0\}\cap \Sigma|^{\frac{1-a}{n+1}}$$
for all $(x,y)\in B_{\eps/2}(\eps k)$.
Finally, using the fact that
$$\lim_{\eps\rightarrow 0} |\{w^\eps_\alpha(x,0) = 0\}\cap\Sigma| = \overline \ell(\alpha)|\Sigma| =0$$
for all $\alpha \leq \alpha_0$ we easily deduce the first inequality in (\ref{eq:w0}).
\end{enumerate} \qed

\vspace{10pt}

\subsection{Proof of Proposition \ref{prop:aux}} 
In order to complete the proof of Proposition \ref{prop:aux}, we construct a corrector $\tw^\eps$ which is equal to $1$ on the (n+1)-dimensional balls $B^+_{r(k,\omega)a^\eps}(\eps k)$.
More precisely, we recall that $D$ is a bounded subset of $\RR^{n+1}_+$, and we introduce
$$\wTe = D\cap \bigcup_{k\in\ZZ^n\cap\Sigma} B^+_{r(k,\omega)a^\eps}(\eps k)$$
and 
$$\wSe = \Sigma\setminus \bigcup_{k\in\ZZ^n\cap\Sigma}  B^n_{r(k,\omega)a^\eps}(\eps k) .$$
We then define a corrector  $\tw^\eps(x,y, \omega)$ which will satisfy all the conditions of Proposition \ref{prop:1}, with the set $\wTe$ instead of $T_\eps$. 
In particular, we will prove that $\tw^\eps$ behaves like $h^\eps_k$ near the $B^+_{r(k,\omega)a^\eps}(\eps k)$.

We consider the following obstacle problem:
\begin{equation}\label{eq:wobst}
\left\{
 \begin{array}{ll} 
\div(y^a\na w )\leq 0 & \mbox{ for }(x,y) \in\RR^{n+1}_+\setminus \wTe \\[5pt]
\ds \lim_{y\rightarrow 0}y^a\pa_y w (x,y) \leq \alpha_0& \mbox{ for } x\in\wSe \\[5pt]
w(x,y) \geq 1&  \mbox{ for } (x,y)\in \wTe \\[5pt]
w (x,0)= 0&  \mbox{ for }x\in  \RR^n\setminus \Sigma,
\end{array}
\right.
\end{equation}
and we define:
$$
 \tw^\eps (x,y,\omega) = \inf\left\{w(x,y,\omega)\,;\, \mbox{$w$ solution of (\ref{eq:wobst})} \right\}.
$$
It is readily seen that $\tw^\eps$ satisfies (\ref{H1bis}).
So in order to complete the proof of Proposition \ref{prop:aux}, we only have to show that 
$\tw^\eps$ is bounded uniformly in $L^\infty(D)$ and that
$\tw^\eps \longrightarrow 0$ in $\Wal$-weak as $\eps$ goes to zero.
\vspace{10pt}

\noindent {\bf Strong convergence in $L^2 (D,|y|^a)$:}\\
First of all, since $\tw^\eps=1=h_{\alpha,k}^\eps (x,y)+o(1)$ on $\wTe$,  (\ref{eq:w0}) implies
$$
w^\eps_{0}(x,y) - o(1)\leq  \tw^\eps  (x,y,\omega)\; \leq\; w^\eps_{0}(x,y) +o(1)\quad  \mbox{ in }  D \quad \mbox{ a.e. } \omega\in\Omega,
$$
which in turn implies (using Proposition  \ref{lem:woe}  again):
\begin{equation}\label{owh}
h_{\alpha,k}^\eps (x,y) -o(1)\; \leq\;  \tw^\eps  (x,\omega)\; \leq\; h_{\alpha,k}^\eps (x,y) +o(1)\quad  \forall (x,y)\in B^+_{\eps/2}(\eps k).
\end{equation}
In particular, we get:
$$||\tw^\eps||_{L^\infty}(\Rn)\leq C.$$

Moreover, a simple computation shows that
$$
\int_{B_\eps(\eps k)\setminus B_{a^\eps}(\eps k)} y^a | h^\eps_{\alpha,k} |^2\, dx \, dy\leq  C \eps^{n+1} 
$$
and it is readily seen that (\ref{owh}) implies
$$
|w_0^\eps(x,y)| \leq C\eps^{1-a} +o(1)=o(1) \qquad \forall (x,y)\in \bigcup_{k\in\ZZ^n} \pa B_{\eps/2}(\eps k) .
$$
We deduce:
$$
||\tw^\eps ||_{L^2(D,|y|^a)} ^2 \leq \sum_{k\in \{\ZZ^n \cap \eps^{-1}\Sigma\}} \int_{B_\eps\setminus B_{a^\eps}} y^a | h^\eps_{\alpha,k} |^2\, dx \, dy + o(1)\int_D |y|^a\, dx\, dy
$$ 
and since $\#\{\ZZ^n \cap \eps^{-1}\Sigma\}\leq C\eps^{-n}$ for all $n$, we have:
\begin{equation}\label{eq:Lp2}
|| \tw^\eps||_{L^2(D,|y|^a)} ^2\leq \eps+o(1)=o(1). 
\end{equation}

In particular
$$
\tw^\eps \longrightarrow 0 \quad \mbox{  in } L^2 (D,|y|^a)-\mbox{strong}.
$$
as $\eps$ goes to zero.

\vspace{10pt}

\noindent {\bf Bound in $\Wa$:}\\ 
Using the definition if $\tw^\eps$ and an integration  by parts, we get:
\begin{eqnarray*}
 \int_{ \Rn \setminus \wTe }y^a |\na \tw^\eps|^2 \, dx \, dy&=&\int_{ \Rn \setminus \wTe }y^a \na \tw^\eps\cdot \na (\tw^\eps-1) \, dx \, dy
\\
& =&  - \int_{\pa \wTe \cup \wSe} [\lim_{y\rightarrow 0} y^a \tw^\eps_y(x,y)] (\tw^\eps(x,y)-1)\, d\sigma(x,y)\\
 & = & - \alpha_0 \int_\wSe  (\tw^\eps(x,0)-1)\, dx
 \end{eqnarray*}
 The $L^\infty$ bound thus yieds
\begin{eqnarray*}
 \int_{ \Rn \setminus \wTe }y^a |\na \tw^\eps|^2 \, dx \, dy&\leq & C \alpha_0 |\wSe| (||\tw^\eps||_{L^\infty} + 1)\leq C,
 \end{eqnarray*}
which completes the proof.
\qed

\vspace{20pt}

\section{Proof of proposition \ref{prop:1}}\label{sec:second}
This section is devoted to the proof of the main proposition. 
We recall that the sets $S_\eps (k,\omega)$ are subsets of  $\RR^n$ with unspecified shapes and  they satisfy
$$
\capa_s (S_\eps(k,\omega)) = \eps^n \gamma(k,\omega).
$$
Lemma \ref{lem:cap} gives
the existence of a function $\vphi^\eps_{k}(x,y,\omega)$ such that
$$
\left\{ \!\! \begin{array}{ll} 
\div(y^a \na \vphi) = 0   & \mbox{ for } (x,y) \in \RR^{n+1}_+\\[5pt]
\vphi(x,0) = 1\quad \qquad\qquad\qquad& \mbox{ for } x \in S_{\eps}(k,\omega) \\[5pt]
\ds \lim_{y\rightarrow 0} y^a \pa \vphi(x,y) = 0 & \mbox{ for } x \notin S_{\eps}(k,\omega)
\end{array}
\right.
$$
and we let $\alpha_0$  and $\tw^\eps (x,y,\omega)$ be given by Proposition \ref{prop:aux}.

We then have:
\begin{enumerate}
\item For a given  $\delta>0$, Lemma \ref{lem:cap} implies that for every $k\in\ZZ^n$ and $\omega\in\Omega$ there exists a constant $R_\delta(k,\omega)$ such that
\begin{equation}\label{eq:gk}
\left| \vphi^\eps_k  (x,\omega) - h^\eps_k(x,y,\omega)  \right|  \leq \delta \, h^\eps_k(x,y,\omega)  \leq \delta\frac{\tgamma(k,\omega)}{R_\delta^{n-1+a}} 
\end{equation}
in $   B^+_{2a^\eps R_\delta }\setminus B^+_{a^\eps R_\delta } (\eps k)$ and for all $\eps>0$.
It is readily seen that for any $R$ there exists $\eps_1(R )$ such that 
\begin{equation}\label{eq:Rd}
a^\eps R  \leq \eps^\sigma /4\quad  \mbox{ for all $\eps\leq \eps_1$.}
\end{equation}
for some $\sigma >1$.

\item Inequality (\ref{eq:hkeps}) in Proposition \ref{prop:aux}  implies that for given $\delta$ and $R$, there exists $\eps_2(\delta,R)<\eps_1(R)$ such that
for all $\eps\leq \eps_2(\delta,R)$, we have
\begin{equation}\label{eq:we}
\left|\tw^\eps (x) - h^\eps_k(x,y,\omega)\right| \leq \frac{\delta}{R^{n-1+a}} \qquad \mbox{ in } B^+_{\eps/2} (\eps k).
\end{equation}
Thanks to (\ref{eq:Rd}), Inequality (\ref{eq:we})  holds in particular in
$B^+_{2a^\eps R }\setminus B^+_{a^\eps R} (\eps k)$.
\end{enumerate}

The  corrector will be constructed by gluing together the functions $\vphi^\eps_k$ (near the sets $S_\eps(k)$) and the function $\tw^\eps $ (away from the sets $S_\eps(k)$).
The gluing has  to be done very carefully so that the corrector satisfies all the properties listed in Proposition   \ref{prop:1}:
For a given $\eps$, we define $\delta_\eps$ to be the smallest positive number such that (\ref{eq:Rd}) and (\ref{eq:we}) hold with $\delta=\delta_\eps$ and $R=R_{\delta_\eps}$.
From the remarks above, we see that $\delta_\eps$ is well defined as soon as $\eps$ is small enough (say smaller than $\eps_2(1,R_1)$).
Moreover, 
for any $\delta>0$, there exists $\eps_0=\eps_2(\delta,R_\delta)$ such that
$\delta_\eps \leq \delta$ for all $ \eps\leq \eps_0.$
In particular 
$$\lim_{\eps\rightarrow 0} \delta_\eps = 0.$$
From now on, we write
$$ R_\eps  = R_{\delta_\eps}.$$

In order to define $w^\eps$, we introduce the cut-off function $\eta_\eps(x,y)$  defined on $D$ and such that
$$ \begin{array}{ll}
\ds \eta_\eps(x,y) = 1 &\ds \mbox{ for  }(x,y)\in D \setminus \bigcup_{k\in\ZZ^n}  B^+_{2 a^\eps R_\eps }(\eps k)\\[10pt]
\ds \eta_\eps(x,y)=0&\ds \mbox{  for }(x,y)\in  \bigcup_{k\in\ZZ^n}   B^+_{a^\eps R_\eps }(\eps k).
\end{array}
$$
We can always choose $\eta$ in such a way that 
$$|\na \eta_\eps|\leq C ( a^\eps R_\eps)^{-1}\quad\mbox{  and }\quad |\Delta \eta_\eps|\leq C(a^\eps R_\eps )^{-2} $$
for $(x,y)\in B^+_{2 a^\eps R_\eps }(\eps k) \setminus   B^+_{a^\eps R_\eps }(\eps k).$
We now set:
$$\cor (x,y)  =   \eta_\eps (x,y) \tw^\eps (x,y)+(1-\eta_\eps(x,y))  \sum_{k\in \ZZ^n\cap D}    \vphi^\eps_{k} (x,y)\,1_{ B^+_{\eps/2}(\eps k)}(x,y) .$$
It satisfies
$$
\cor (x,y,\omega)=\left\{
\begin{array}{ll}
\vphi^\eps_{k} (x,y)& \mbox{ for }(x,y)\in   D\cap B^+_{a^\eps R_\eps }(\eps k) \quad \forall k\in\ZZ^n \\[8pt]
\tw^\eps (x,y) &\ds  \mbox{ for } (x,y)\in D \setminus \bigcup_{k\in\ZZ^n}  B^+_{2 a^\eps R_\eps }(\eps k).
\end{array}
\right.
$$
To simplify the notations in the sequel, we denote
$$ 
\vphi^\eps (x,y):= 
 \sum_{k\in \ZZ^n\cap D} \vphi^\eps_{k} (x,y,\omega)\,1_{B^+_{\eps/2}(\eps k)}(x,y).
$$

The properties of $\cor$ are summarized in the following lemma, which implies Proposition \ref{prop:1}:
\begin{lemma}
The function $\cor$ satisfies the following properties:
\item [(i)] $ \cor(x,0) = 1$ for $x\in S_\eps$ and $||\cor||_{L^\infty(D)}\leq C$.
\item [(ii)] $ \cor$ converges to zero in $L^2(D,|y|^a)$-strong as $\eps$ goes to zero.
\item[(iii)] $ \cor$ is bounded in $\Wa$.
\item  [(iv)] $ \cor$ satisfies  (\ref{H1'}).
\end{lemma}
{\it Proof:}
\begin{enumerate}
\item[(i)] Immediate consequence of the definition of $\cor$ since $\vphi^\eps_{k} = 1$ on $S_\eps(k,\omega)$ and $\tw^\eps$ and $\vphi_k^\eps$ are bounded in $L^\infty$.
\item[(ii)] 
Since 
$S_\eps(k,\omega) \subset B^n_{a^\eps M}(\eps k)$, we have:
$$  \vphi^\eps_{k}(x,y,\omega) \leq C \eps^n \frac{M^{n-1+a}}{\nu_{n-1+a}} h(x-\eps k,y)$$
for all $(x,y)$ such that $|(x-\eps k,y)| \geq a^\eps M $.
Since $\vphi^\eps_{k}\leq 1$ in $B_{a^\eps M }(\eps k)$, we get
\begin{eqnarray*}
&& \int_D |y|^a \left|(1-\eta_\eps) \vphi^\eps  \right|^2\, dx\, dy \\
& &\quad \leq   \sum_{k\in \ZZ^n\cap \eps^{-1} \Sigma} \int_{B_{2a^\eps R}(\eps k)} |y|^a \left|\vphi^\eps_k  \right|^2\, dx\, dy \\
& &\quad  \leq    \sum_{k\in \ZZ^n\cap\eps^{-1}  \Sigma} \int _{B_{a^\eps M}(\eps k)}\!\!\!\!\!\!\!\!\!|y|^a \, dx\,dy  \\
&& \qquad + C \sum_{k\in \ZZ^n\cap\eps^{-1}  \Sigma} \int _{B_{2a^\eps(R_\eps)(\eps k)\setminus B_{a^\eps M}(\eps k)}}\!\!\!\!\!\!\!\!\!\!\!\!  |y|^a\left( \eps^n    \frac{M^{n-1+a}}{\nu_{n-1+a}}           h(x-\eps k)\right)^2 \, dx \\
& & \quad\leq   \sum_{k\in \ZZ^n\cap\eps^{-1}  \Sigma} (a^\eps M)^{n+1+a}  \\
&& \qquad + C \sum_{k\in \ZZ^n\cap\eps^{-1}  \Sigma}  \eps^{2n}  (a^\eps M )^{n+1-2(n-1+a)} M^{n-1+a}
\end{eqnarray*}
Using (\ref{eq:Rd}) and the definition of $a^\eps$, we deduce:
\begin{eqnarray*}
\left\| (1-\eta_\eps) \vphi^\eps \right\|^2_{L^2(D,|y^a|)} 
&\leq & C  (M) \eps^\frac{2n-an }{n-1+a} .
\end{eqnarray*}
Estimate (\ref{eq:Lp2}) thus implies
\begin{eqnarray*}
|| \cor ||_{L^2(D,|y^a|)}  & \leq &  || \tw^\eps  ||_{L^2(D,|y^a|)} + || (1-\eta_\eps) \vphi^\eps  ||_{L^2(D,|y^a|)}   = o(1).
\end{eqnarray*}
and therefore
$$\cor \longrightarrow 0 \qquad L^2(D,|y|^a) \mbox{-strong.}$$

\item[(iii)] Next, we want to show that $\cor$ is bounded in $\Wa$. 
First, we note that outside $\cup_{k\in\ZZ^n}B_{\eps/2}(\eps k)$ we have $\na \cor =\na \tw^\eps$ which is bounded in $\Wa$.
Next, we see that 
in $B_{\eps/2}(\eps k)$, we have:
\begin{equation}\label{eq:grad}
\nabla \cor = \na \eta_\eps (\tw^\eps -\vphi^\eps_{k}) +  \eta_\eps \na  \tw^\eps+(1-\eta_\eps) \na \vphi^\eps_{k} 
\end{equation}
Since $\tw^\eps$ and  $\vphi^\eps $ are both bounded in $\Wa$ (thanks to (\ref{eq:gradg})), we see that in order to show that $\nabla \cor $ is bounded in $L^2(D,|y|^a)$, we only have to show that
$$
\int_D y^a | \na \eta_\eps (\tw^\eps -\vphi^\eps) |^2 \, dx\, dy \leq C.
$$

For that purpose, we notice that (\ref{eq:gk}) and (\ref{eq:we}) yield
$$|\tw^\eps -\vphi^\eps_{k} | \leq C\frac{\delta_\eps}{R_\eps^{n-1+a}}\quad  \mbox{ in } B_{2R_\eps a^\eps }(\eps k)\setminus B_{R_\eps a^\eps}(\eps k), $$ 
and so,  using the definition of $\eta_\eps(x,y)$, we deduce:
\begin{eqnarray*}
&& \int_D y^a | \na \eta_\eps (\tw^\eps -\vphi^\eps) |^2 \, dx \, dy\\
&&  \qquad\qquad \leq  
\sum_{k \in \eps\ZZ^n\cap \Sigma} \int_{B_{2 R_\eps a^\eps}(\eps k)\setminus B_{ R_\eps a^\eps}(\eps k)} y^a | \na \eta_\eps (\tw^\eps -\vphi^\eps_{k}) |^2 \, dx \\
&&\qquad\qquad  \leq  
\sum_{k \in \eps\ZZ^n\cap \Sigma} (R_\eps a^\eps)^{n+1+a}(R_\eps a^\eps)^{-2} \frac{\delta_\eps ^2}{R_\eps ^{2(n-1+a)}} \\
&&\qquad\qquad  \leq  
\sum_{k \in \eps\ZZ^n\cap \Sigma} R_\eps ^{-(n-1+a)} \eps^n \delta_\eps ^2 \\
&& \qquad\qquad\leq C\eps^{-n} \eps^n \delta_\eps = C \delta_\eps,
\end{eqnarray*}
where we used the fact that  we can always assume that $\delta_\eps<1$ and $R_\eps\geq 1$.
For latter use, we note that we actually proved
\begin{equation}\label{eq:conv2}
 \int_D y^a | \na \eta_\eps (\tw^\eps -\vphi^\eps) |^2 \, dx \, dy\longrightarrow 0 \mbox{ when }\eps\rightarrow 0.
\end{equation}
\vspace{10pt}

\item[(iv)] It remains to show that (\ref{H1'}) holds. We only show the inequality (the equality follows  easily).
Let $v^\eps$ be a sequence of functions satisfying:
$$
\left\{ 
\begin{array}{l}
v^\eps (x,0)\geq 0\quad \mbox{ for } x\in T_\eps\\[3pt]
||v^\eps||_{L^\infty(D)} \leq C\\[3pt]
v^\eps\longrightarrow v\quad \mbox{ in } \Wa-\mbox{weak.}
\end{array}
\right.
$$
Then for any $\phi\in\mathcal D(D)$, we have:
\begin{eqnarray*}
&& \ds  -\int_D y^a \na  w^\eps \cdot \na v^\eps \phi\, dx\, dy \\
& & \quad= -\int_D y^a \na \eta _\eps \cdot \na v^\eps (\tw^\eps-\vphi^\eps )\phi\, dx\, dy -\int_D y^a \na \tw^\eps \cdot \na v^\eps \phi \,\eta_\eps\, dx\, dy \\
&&\qquad -\int_D y^a \na \vphi^\eps \cdot \na v^\eps\phi \, (1-\eta_\eps) \, dx\, dy \\
& & \quad= -\int_D y^a \na \eta _\eps \cdot \na v^\eps (\tw^\eps-\vphi^\eps )\phi\, dx\, dy \\
&&\qquad + \int_\Sigma  (\lim_{y\rightarrow 0} y^a\pa_y \tw^\eps) v^\eps \phi\, \eta_\eps\, dx+\int_\Sigma ( \lim_{y\rightarrow 0} y^a\pa_y \vphi^\eps) v^\eps \phi\, (1-\eta_\eps)\, dx\\
&& \qquad+ \int_D y^a \na \tw^\eps \cdot  \na( \phi \eta_\eps)v^\eps \, dx\, dy 
 +\int_D y^a \na \vphi^\eps \cdot  \na(\phi  (1-\eta_\eps))v^\eps \, dx\, dy
 \end{eqnarray*}
where we used the fact that $\div(y^a \na \tw^\eps) = 0$ on $\mbox{supp } \eta^\eps$ and 
$\div(y^a \na \vphi^\eps) = 0$ on $\mbox{supp } (1-\eta^\eps)$. 
The first term goes to zero thanks to (\ref{eq:conv2}) and the weak convergence of $\na v^\eps$ in $L^2(D,|y|^a)$, and the boundary terms satisfy
\begin{eqnarray*}
\lim_{\eps\rightarrow 0} \int_\Sigma  (\lim_{y\rightarrow 0} y^a\pa_y \tw^\eps) v^\eps \phi \eta_\eps\, dx
&=  &\lim_{\eps\rightarrow 0}\int_\Sigma  \alpha_0 v^\eps \phi \eta_\eps\, dx\\
&= &\lim_{\eps\rightarrow 0} \int_\Sigma  \alpha_0 v \phi \eta\, dx\\
\end{eqnarray*}
and
\begin{eqnarray*}
\lim_{\eps\rightarrow 0}\int_\Sigma  (\lim_{y\rightarrow 0} y^a\pa_y \vphi^\eps) v^\eps \phi (1-\eta_\eps)\, dx
& =& \lim_{\eps\rightarrow 0}\int_{T_\eps}  (\lim_{y\rightarrow 0} y^a\pa_y \vphi^\eps) v^\eps \phi (1-\eta_\eps)\, dx\\
&\leq& 0.
\end{eqnarray*}

Finally, the last two terms can be rewritten as:
\begin{eqnarray*}
& &   \int_D y^a \na \tw^\eps \cdot  \na( \phi \eta_\eps)v^\eps\, dx\, dy  +\int_D y^a \na \vphi^\eps \cdot \na(\phi  (1-\eta_\eps))  v^\eps \, dx\, dy  \\
& & \quad = \int_D y^a \na ( \tw^\eps-\vphi^\eps)  \cdot (\na  \eta_\eps)\,  v^\eps \phi \, dx\, dy \\
&&   \qquad + \int_D y^a \, v^\eps \eta_\eps \, \na \tw^\eps \cdot  \na \phi \,  dx\, dy  +\int_D y^a \, v^\eps (1-\eta_\eps) \,  \na \vphi^\eps \cdot  \na \phi \, dx\, dy 
\end{eqnarray*}
Using the weak convergence of $\na \tw^\eps $ and $\na \vphi^\eps$ to zero, we see that in order to prove (\ref{H1'}), it only remains to prove that
$$\int_D y^a \na ( \tw^\eps-\vphi^\eps)  \cdot (\na  \eta_\eps)\,  v^\eps \phi \, dx\, dy \longrightarrow 0 \mbox{ when }\eps\rightarrow 0.
$$
Since $v^\eps$ is bounded in $L^\infty$, it is enough to show that
$$\int_D |y|^a |\na ( \tw^\eps-\vphi^\eps)|  \, |\na  \eta_\eps| \, dx\, dy \longrightarrow 0 \mbox{ when }\eps\rightarrow 0.
$$

For that purpose, we recall that 
$$|\tw^\eps  - \vphi^\eps_{k} | \leq \frac{\delta_\eps }{R_\eps ^{n-1+a}}\qquad  \mbox{ in }B^+_{2a^\eps R_\eps }\setminus B^+_{a^\eps R_\eps } ,$$
and
$$
\left\{
\begin{array}{l}
\div(y^a \na ( \tw^\eps  - \vphi^\eps_{k} )) = 0 \qquad \qquad \mbox{ for }(x,y)\in  B^+_{4a^\eps R }\setminus  B^+_{a^\eps R_\eps /2} \\[8pt]
\ds \lim_{y\rightarrow 0} y^a \pa_y  ( \tw^\eps  - \vphi^\eps_{k} ) (x,y)= \alpha_0  \qquad \mbox{ for }x \in  B^n_{4a^\eps R }\setminus  B^n_{a^\eps R_\eps /2}.
\end{array}
\right.
$$
In particular, interior gradient estimates  (see \cite{CSS}) implies
$$  | \na(\tw^\eps  - \vphi^\eps_{k})| \leq \frac{\delta_\eps }{R_\eps ^{n-1+a}} (a^\eps R_\eps)^{-1}  + C (a^\eps R_\eps)^{-a}$$
in $B^+_{2a^\eps R_\eps } \setminus B^+_{a^\eps R_\eps} $.
We deduce:
\begin{eqnarray*}
&& \int_D |y|^a |\na ( \tw^\eps-\vphi^\eps)|  \, |\na  \eta_\eps| \, dx\, dy  \\
&&\quad \leq\sum_{k \in \eps\ZZ^n\cap \Sigma}   \int_{B^+_{2a^\eps R_\eps } \setminus B^+_{a^\eps R_\eps} }  |y|^a |\na ( \tw^\eps-\vphi^\eps)|  \, |\na  \eta_\eps| \, dx\, dy\\
&&\quad \leq\sum_{k \in \eps\ZZ^n\cap \Sigma}    \frac{C \delta_\eps }{R_\eps ^{n-1+a}}  (a^\eps R_\eps )^{-2}  \int_{B^+_{2a^\eps R_\eps } \setminus B^+_{a^\eps R_\eps} }  |y|^a \, dx\, dy\\
&&\qquad +\sum_{k \in \eps\ZZ^n\cap \Sigma}   C(a^\eps R_\eps )^{-1-a}  \int_{B^+_{2a^\eps R_\eps } \setminus B^+_{a^\eps R_\eps} }  |y|^a \, dx\, dy\\
&&\quad \leq\sum_{k \in \eps\ZZ^n\cap \Sigma}    \frac{\delta_\eps }{R_\eps ^{n-1+a}}  (a^\eps R_\eps )^{-2}  (a^\eps R_\eps )^{n+1+a}  \\
&&\qquad +\sum_{k \in \eps\ZZ^n\cap \Sigma}  C(a^\eps R_\eps )^{-1-a}  (a^\eps R_\eps )^{n+1+a}  \\
&&\quad \leq  \frac{C \delta_\eps }{R_\eps ^{n-1+a}} \eps^{-n}   (a^\eps R_\eps )^{n-1+a}  
+ C \eps^{-n}   (a^\eps R_\eps )^{n}  .
\end{eqnarray*}
Using (\ref{eq:Rd}) and the definition of $a^\eps$, we deduce:
$$ \int_D |y|^a |\na ( \tw^\eps-\vphi^\eps)|  \, |\na  \eta_\eps| \, dx\, dy  
\leq  C \delta_\eps  + C \eps^{n(\sigma-1)}.
$$
which concludes the proof since $\sigma>1$ and $\lim_{\eps\rightarrow 0}\delta_\eps = 0$.
\end{enumerate}
\qed

\vspace{10pt}

\noindent{\bf Acknoledgment:} L. Caffarelli was partially supported by NSF grant DMS-0140338. A. Mellet was partially supported by NSERC discovery grant 315596.

\newpage

\appendix

\section{Proof of Lemma \ref{lem:cap}}\label{app:capa}
We now turn to the proof of Lemma  \ref{lem:cap}. We take $k=0$ and we recall that $\vphi^\eps_0$ is the capacity potential associated to $S_\eps(0)$. It satisfies (\ref{eq:vphik}) and (\ref{eq:gradg}).
We then introduce the function
$$ G(x,\xi,y,\tau) = h(x-\xi,y-\tau)+h(x-\xi,y+\tau)$$
which satisfies
$$ \div_{\xi,\tau}(|\tau|^a\na_{\xi,\tau} G) =- \mu_{n,a}\delta{(x-\xi,y-\tau)}-\mu_{n,a}\delta{(x-\xi,y+\tau)}$$
and
$$ \lim_{\tau\rightarrow0} \tau^a \pa_\tau G (x,\xi,y,\tau)=0$$
for all $x$, $\xi$ and $y$.
If $y>0$, we deduce that for any function $\vphi(x,y)$, we have:
\begin{eqnarray*}
&&\int_{\tau>0} \tau^a \na_{\xi,\tau} G(x,\xi,y,\tau) \na_{\xi,\tau} \vphi(\xi,\tau)\, d\xi\, d\tau.\\
& &\qquad\qquad = - \int_{\tau>0}\div( \tau^a \na_{\xi,\tau} G(x,\xi,y,\tau))\vphi(\xi,\tau)\, d\xi\, d\tau \\
&&\qquad\qquad\qquad-\int_{\RR^n}  \lim_{\tau\rightarrow0} \tau^a \pa_\tau G (x,\xi,y,\tau) \vphi(\xi,0)\, d\xi\\
& &\qquad\qquad = \mu_{n,a}\vphi(x,y).
\end{eqnarray*}
Moreover, if $\vphi_0^\eps(x,y)$ is the capacity potential associated to $S_\eps(0)$, then (\ref{eq:vphik}) yields
\begin{eqnarray*}
&&\int_{\tau>0} \tau^a \na_{\xi,\tau} G(x,\xi,y,\tau) \na_{\xi,\tau} \vphi_0^\eps(\xi,\tau)\, d\xi\, d\tau.\\
& &\qquad\qquad =
- \int_{\tau>0} G(x,\xi,y,\tau)\div( \tau^a \na_{\xi,\tau} \vphi_0^\eps(\xi,\tau))\, d\xi\, d\tau \\
& &\qquad\qquad\qquad -\int_{\RR^n}G (x,\xi,y,0)  \lim_{\tau\rightarrow0} \tau^a \pa_\tau \vphi_0^\eps(\xi,\tau)  \, d\xi\\
& &\qquad\qquad= -2 \int_{\RR^n}  h(x-\xi,y) \lim_{\tau\rightarrow0} \tau^a \pa_\tau \vphi_0^\eps(\xi,\tau) \, d\xi.
\end{eqnarray*}
Combining those two equalities, we get:
$$ \mu_{n,a}\vphi_k^\eps(x,y)=-2 \int_{S_\eps(0)}  h(x-\xi,y) \lim_{\tau\rightarrow0} \tau^a \pa_\tau \vphi_0^\eps(\xi,\tau) \, d\xi.$$
Next, we note that (\ref{eq:gradg}) yields,  after integration by parts and using (\ref{eq:vphik}):
$$ \eps^n \gamma (0)=  \int_{\RR^n} \tau^a |\na \vphi_0^\eps(\xi,\tau)|^2 \, d\xi=- \int_{S_\eps(0)}   \lim_{\tau\rightarrow0} \tau^a \pa_\tau \vphi_0^\eps (\xi,\tau) \, d\xi,$$
and therefore
\begin{eqnarray*}
&& \vphi_0^\eps(x,y)-\frac{2}{\mu_{n,a}}\eps^n\gamma(0) h(x,y)\\
&&\qquad\qquad=-\frac{2}{\mu_{n,a}} \int_{S_\eps(0)}  [h(x-\xi,y) -h(x,y)]\lim_{\tau\rightarrow0} \tau^a \pa_\tau \vphi_0^\eps(\xi,\tau) \, d\xi.
\end{eqnarray*}

In order to conclude, we recall that $S_\eps(0) \subset B_{M a^\eps}(0)$ and so we have $|\xi|\leq M a^\eps$  in the previous integral. If $(x,y)$ is such that $|(x,y)|\geq R a^\eps$ with $R\geq 8 M$, we deduce that for all $\xi\in S_\eps(0)$, we have:
\begin{eqnarray*}
|h(x-\xi,y) -h(x,y)| & \leq& \sup_{\xi^*\in B_{M a^\eps}(0)} |\na_{x,y} h(x-\xi^*,y)| |\xi| \\
& \leq& \sup_{\xi^*\in B_{M a^\eps}(0)} \frac{|\xi|}{((x-\xi^*)^2+y^2)^{\frac{n-a}{2}}} \\
& \leq&  \frac{C |\xi|}{(x^2+y^2)^{\frac{n-a}{2}}} \\
& \leq& \frac{C |\xi|}{(x^2+y^2)^{\frac{1}{2}}} h(x,y) \\
& \leq& \frac{C M}{R} h(x,y) 
\end{eqnarray*}
We can thus write
\begin{eqnarray*}
&&\left| \vphi_0^\eps(x,y)-\frac{2}{\mu_{n,a}}\eps^n\gamma(0) h(x,y)\right| \\
&&\qquad\qquad\qquad\leq   \frac{C M}{R}\frac{2}{\mu_{n,a}} h(x,y) \int_{S_\eps(0)} \left|  \lim_{\tau\rightarrow0} \tau^a \pa_\tau \vphi_0^\eps(\xi,\tau)\right| \, d\xi\\
&&\qquad\qquad\qquad\leq  \frac{C M}{R}\frac{2}{\mu_{n,a}} \eps^n \gamma (0) h(x,y)
\end{eqnarray*}
where the right hand side is bounded by $ \delta \frac{2}{\mu_{n,a}}  \eps^n \gamma (0) h(x,y)
$ if $R$ is large enough.
\qed

\section{Proof of Lemma \ref{lem:alpha}.} \label{sec:alpha}
\noindent (i) 
For a given set $A$, it is readily seen from the definition of $\ov_{\alpha,A}$ that if $  \alpha' \leq \alpha$, then $\ov_{\alpha',A}$ is admissible for the obstacle problem with $\alpha$: It follows that
$$ \ov_{\alpha,A} \leq \ov_{\alpha',A}\qquad  \mbox{ for any $\alpha$, $\alpha'$ such that } \alpha' \leq 
\alpha$$
and so $\alpha\mapsto m_\alpha(A,\omega)$ is nondecreasing. The result follows from the definition  of $\oell(\alpha)$.
\vspace{10pt}

\noindent (ii)
If $\alpha$ is negative, then we have
$$\lim_{y\rightarrow 0} y^a \pa_y \ov_{\alpha,tB} (x,y) <0\quad\mbox{ for } x\in\RR^n .$$
Since
$\ov_{\alpha,tB} (x,y)\geq 0$ for $(x,y)\in\RR^{n+1}_+,$
we deduce
$$ \ov_{\alpha,tB} (x,0) >0 \quad\mbox{ for } x\in\RR^n.$$ 
It follows that $m_\alpha (tB,\omega) =0$  for all $t>0$, and so $\oell(\alpha)=0$ for all $\alpha<0$.
\vspace{10pt}

If $r(k,\omega)$ is bounded below:
$$\;\; r(k,\omega) \geq \underline r >0 \mbox{ for all $k\in\ZZ^n$, a.e. $\omega\in\Omega$}, $$
then, 
we define
$$ \vphi(x,y) =\frac{\underline r^{n-1+a} }{(|x|^2+y^2)^\frac{n-1+a}{2}} -  \alpha \int_{B^n_1(0)} \frac{\nu_{n+1+a}}{(|x-x'|^2+y^2)^{\frac{n-1+a}{2}}}\, dx' -C_0$$
with 
$$C_0 =\underline r^{n-1+a}- \alpha \int _{B^n_1(e)} \frac{\nu_{n+1+a}}{|z|^{n-1+a}}\, dz
$$
where $e $ denote any unit vector in $\RR^n$.
In particular, we have
$$ \vphi(x,0) = 0 \mbox { if } |x|=1, $$
and, if $\alpha$ is small enough
$$
\vphi(x,0)> 0 \mbox { if } |x|<1, 
$$
and 
$$
\vphi(x,y)< 0 \mbox { if } |x|=1,\quad y>0
$$
(we note that $\vphi$ is the sum of a term which is decreasing with respect to $|x|$ and one which is increasing).
Since $\vphi$ satisfies
$$\lim_{y\rightarrow 0}y^a \pa_y \vphi(x,y) = \alpha - \underline \gamma \delta(x) \geq \alpha - \tgamma(0,\omega) \delta(x),$$
for all $x\in B^n_1(0)$,
we deduce
$$\ov_{\alpha,tB} (x,0)> \vphi(x,0) >0 \qquad \mbox{ in } B^n_1(0).$$
Since we can do this in any ball $B^n_1(k)$, we must have
 $m_\alpha (tB^n,\omega) =0$ for all $t>0$, and so $\oell(\alpha)=0$ for all $\alpha$ small enough.
\vspace{10pt}

\noindent (iii) 
We consider the function
$$\psi(x,y) = \frac{\overline r^{n-1+a} }{(|x|^2+y^2)^\frac{n-1+a}{2}} -\alpha \int_{B^n_1(0)} \frac{\nu_{n+1+a}}{(|x-x'|^2+y^2)^{\frac{n-1+a}{2}}}\, dx' +C,
$$
where the constant $C$ will be chosen later.
It satisfies 
$$\lim_{y\rightarrow 0}y^a\pa_y \psi(x,y)= \alpha - \overline \gamma\, \delta(x) \leq \alpha - \tgamma(0,\omega) \delta(x)\qquad \forall x\in B^n_1(0),$$
$$ \psi(x,y)\longrightarrow C \qquad \mbox{when  } |x|^2+y^2\rightarrow \infty.$$
and we note that $\psi(x,0)$ is radially symmetric. Moreover, when $\alpha$ is such that
$$ 
\alpha  \int_{B_1(0)} \frac{\nu_{n+1+a}}{|e_1-x'|^{n-1+a}}\, dx'
\geq \overline r^{n-1+a} 
$$
then
$$\psi(x,0)<C\qquad \mbox{ when } |x|=1.$$
Since $\div(y^a\na \psi)=0$ for $y>0$ and $\lim_{y\rightarrow 0}y^a \pa_y \psi(x,y)=0$ for $x\notin B^n_1(0)$,  the strong maximum principle and Hopf Lemma yield that the minimum of $\psi(x,y)$ is reached for $y=0$ and $x\in B^n_1(0)$, and with an appropriate choice of the constant $C$, we can always assume that this minimum is $0$:
$$\inf_{\RR^{n+1}_+} \psi(x,y) = \inf_{B^n_1(0)} \psi(x,0)=0$$
Finally,
if  $\alpha$ is such that
$$\alpha\int_{B^n_1(0)}\left[ \frac{1}{|\frac{e_1}{2}-x'|^{n-1+a}}- \frac{1}{|e_1-x'|^{n-1+a}} \right] dx' \geq \overline r^{n-1+a} (4^{n+1}-1)
$$
then $\psi (x,0)$ reaches its minimum when $|x|=R_\alpha $ with $R_\alpha<1/4$.

We now consider the function $\vphi(x,y)$ defined by:
$$ \vphi(x,y) = 
\left\{
\begin{array}{ll}
\psi(x-k,y) & \mbox{ for }(x,y)\in B^+_{1/4}(k) \\[5pt]
\inf_{k'}\psi(x-k',y) & \mbox{ for  }(x,y)\in \RR^{n+1}_+\setminus \cup_{k'}  B^+_{1/4}(k') 
\end{array}
\right.
$$
We clearly have
$$\lim_{y\rightarrow 0}y^a \pa_y \vphi(x,y) \leq \alpha_0 -\gamma(k,\omega) \delta(x-k) \qquad \mbox{ for }x\in B^n_{1/4}(k) $$
and 
$$\lim_{y\rightarrow 0}y^a \pa_y \vphi(x,y) \leq \alpha_0 \qquad  \mbox{ for }x\in \RR^{n}\setminus \cup_{k'} B^n_{1/4}(k') .$$
In order to prove that $\vphi$ is a  supersolution for the obstacle problem, we only have to check that 
$$\psi(x-k,y) = \inf_{k'}\psi(x-k',y) \quad \mbox{ for }(x,y)\in \pa B^+_{1/4}(k)$$
or equivalently
$$\psi(x,y) = \inf_{k'}\psi(x-k',y)\quad  \mbox{ for }(x,y)\in\pa B^+_{1/4}(0).$$
It is readily seen that this amounts to showing that
$$ \alpha \, \nu_{n+1+a} \int_{B^n_1(0)} \frac{1}{(|x-x'|^2+y^2)^{\frac{n-1+a}{2}}}-\frac{1}{(|x-k-x'|^2+y^2)^{\frac{n-1+a}{2}}}\, dx' \geq 4^{n-1}\overline r^{n-1+a} $$
for all $k\in \ZZ^n\setminus \{0\}$ and all $(x,y)\in \pa B^+_{1/4}(0)$. 
This inequality is obviously satisfied  if $\alpha$ is large enough provided we can prove that
$$
 \int_{B^n_1(0)} \frac{1}{(|x-x'|^2+y^2)^{\frac{n-1+a}{2}}}-\frac{1}{(|x-k-x'|^2+y^2)^{\frac{n-1+a}{2}}}\, dx' >0 $$
for all $k\in \ZZ^n\setminus \{0\}$ and all $(x,y)\in \pa B^+_{1/4}(0)$. 
This is equivalent to 
$$
 \int_{B^n_1(x)} \frac{1}{(|x'|^2+y^2)^{\frac{n-1+a}{2}}}\, dx'>  \int_{B^n_1(x-k)}  \frac{1}{(|x'|^2+y^2)^{\frac{n-1+a}{2}}}\, dx' >0 $$
whick holds for all  $(x,y)\in \pa  B^+_{1/4}(0)$ since $|k|\geq 1$.

By definition of $\ov_{\alpha,tB}$, we deduce that
$$ \ov_{\alpha,tB}(x) \leq \vphi(x,y)\quad \mbox{ in } t B^n \mbox{ a.s.} $$
In particular, this implies that $ \ov_{\alpha,tB^n}$ vanishes in $t B^n \setminus \cup_{k\in\ZZ^n} B_{1/2}(k)$,
and so 
$$ \frac{\overline m_\alpha(t B^n,\omega) }{|tB^n|}\geq  \left(\frac{|C_1| -|B^n_{1/2}|}{|C_1|}\right)  = 1-\frac{\omega_n}{2^n}\quad\mbox{ a.s. }$$  
We conclude
$$ \oell(\alpha ) \geq 1-\frac{\omega_n}{2^n}>0.
$$ \qed

\bibliography{bibli}

\end{document}